\documentclass[11 pt]{amsart}
          
\usepackage{fullpage,eucal,amsmath,amsthm,amsfonts,verbatim,amsopn,amssymb}

\usepackage[all]{xy}

\newcommand\barr{\begin{array}}
\newcommand\earr{\end{array}}
\newcommand\bit{\begin{itemize}}
\newcommand\eit{\end{itemize}}
\newcommand\bce{\begin{center}}
\newcommand\ece{\end{center}}

\newcommand{\ignore}[1]{}

\newtheorem{theorem}{Theorem}[section]
\newtheorem{conjecture}[theorem]{Conjecture}
\newtheorem{lemma}[theorem]{Lemma}
\newtheorem{corollary}[theorem]{Corollary}
\newtheorem{proposition}[theorem]{Proposition}

\theoremstyle{definition}

\newtheorem{definition}[theorem]{Definition}
\newtheorem{example}[theorem]{Example}
\newtheorem{remark}[theorem]{Remark}
\newtheorem{remarks}[theorem]{Remarks}

\newtheorem{question}[theorem]{Question}

\newcommand\bthm{\begin{theorem}}
\newcommand\ethm{\end{theorem}}
\newcommand\bcn{\begin{conjecture}}
\newcommand\ecn{\end{conjecture}}
\newcommand\bla{\begin{lemma}}
\newcommand\ela{\end{lemma}}
\newcommand\bco{\begin{corollary}}
\newcommand\eco{\end{corollary}}
\newcommand\bpro{\begin{proposition}}
\newcommand\epro{\end{proposition}}
\newcommand\bdf{\begin{definition}}
\newcommand\edf{\end{definition}}
\newcommand\bex{\begin{example}}
\newcommand\eex{\end{example}}
\newcommand\brm{\begin{remark}}
\newcommand\erm{\end{remark}}
\newcommand\brms{\begin{remarks}}
\newcommand\erms{\end{remarks}}
\newcommand\bqu{\begin{question}}
\newcommand\equ{\end{question}}

\newcommand\bprf{\begin{proof}}
\newcommand\eprf{\end{proof}}

\newcommand\msk{\medskip}

\newcommand\sub{\subseteq}
\renewcommand\sup{\supseteq}

\newcommand\ra{\rightarrow}

\newcommand\lra{\longrightarrow}
\newcommand\lla{\longleftarrow}

\newcommand\al{{\alpha}}
\newcommand\be{{\beta}}
\newcommand\ga{{\gamma}}
\newcommand\dl{\delta}

\newcommand\cm{\otimes}

\newcommand\Om{\Omega}

\newcommand\pd{\partial}

\newcommand\xbar{\overline{x}}
\newcommand\ybar{\overline{y}}

\newcommand\abar{\overline{a}}
\newcommand\bbar{\overline{b}}
\newcommand\cbar{\overline{c}}
\newcommand\albar{\alpha}
\newcommand\bebar{\beta}
\newcommand\gabar{\gamma}

\newcommand\A{{\mathbb{A}}}

\newcommand\N{{\mathbb{N}}}
\newcommand\Z{{\mathbb{Z}}}

\newcommand\C{{\mathbb{C}}}

\newcommand\DD{{\mathcal{D}}}

\newcommand\calO{{\mathcal{O}}}
\renewcommand\O{{\mathcal{O}}}

\newcommand\Dun{{\underline{D}}}
\newcommand\dun{{\underline{d}}}

\newcommand\Hom{\mathrm{Hom}}
\newcommand\Alg{\mathrm{Alg}}

\newcommand\Spec{\mathrm{Spec}}
\newcommand\BSpec{\mathbf{Spec}}
\newcommand\DSpec{\DD\mbox{-}\mathrm{Spec}}
\newcommand\Der{\mathrm{Der}}
\newcommand\Ker{\mathrm{Ker}}

\newcommand\HS{\mathrm{HS}}

\newcommand\dring{$\DD$-ring}
\newcommand\drings{$\DD$-rings}
\newcommand\dfield{$\DD$-field}

\newcommand\dsch{$\DD$-\textup{scheme}}

\begin{document}

\title{Prolongations in differential algebra}
\author{Eric Rosen}
\address{Department of Mathematics\\
Massachusetts Institute of Technology\\
77 Massachusetts Ave.\\
Cambridge, MA 02139}
\email{rosen@math.mit.edu}
\urladdr{http://math.mit.edu/$\sim$rosen}
\thanks{This paper constitutes part of my dissertation, written at the University 
of Illinois at Chicago.  I would like to thank my advisor, David Marker,
for his guidance and support, and Lawrence Ein and Henri Gillet for 
helpful discussions on this material.  
I am also grateful to the referee for a careful reading of the paper and many
detailed comments that have greatly improved the presentation.}

\date{\today}

\begin{abstract}
We develop the theory of higher prolongations of algebraic
varieties over fields in arbitrary characteristic with commuting Hasse-Schmidt derivations.
Prolongations were introduced by Buium in the context of fields of characteristic 0 with
a single derivation.  Inspired by work of Vojta, we give a new construction of higher 
prolongations in a more general context.
Generalizing a result of Buium in characteristic 0, we prove that these 
prolongations are represented by a certain functor, which shows that they can be viewed
as `twisted jet spaces.'
We give a new proof of a theorem of Moosa, Pillay, and Scanlon
that the prolongation functor and jet space functor commute.  
We also prove
that the $m^{th}$-prolongation and $m^{th}$-jet space of a variety
are differentially isomorphic by showing that their infinite prolongations are 
isomorphic as schemes.
\end{abstract}

\maketitle

\section*{Introduction}

Prolongations of algebraic varieties over a differential field of characteristic 0 were introduced by 
Buium~\cite{Bui92}, 
and have also been considered in more general contexts~\cite{BV95,BV96,Sca97,MPS}.
The purpose of this paper is to develop the basic theory of prolongations of algebraic
varieties over fields 
with finitely many commuting Hasse-Schmidt (or `higher') derivations.  Let us begin by describing the idea 
behind Buium's construction and the connection to jet spaces of varieties.  We then
describe the content of the paper in more detail.

Let $(K,\dl)$ be a differential field and let $R = K\{x_1, \ldots, x_m\}$ be the ring of differential
polynomials, which is the polynomial ring in the infinitely many variables 
$\dl^n x_i, 0 \leq n$ and $i \leq m$.  We say that $f \in R$ has order
$\leq n$ if for every variable $\dl^j x_i$ that occurs in $f$, one has $j \leq n$.  
Observe that the set of elements of order $\leq n$ is a subring of $R$.
Elements of $R$ can be viewed as functions on 
affine $m$-space $\A^m = K^m$ in a natural way.  

More generally, let $X \sub \A^m$ be an affine variety over $(K,\dl)$ and let 
$A $ be the ring of regular functions on $X$.  As above, there is a 
natural way to define the ring of differential polynomial functions of order $\leq n$ on $X$, 
called the $n^{th}$-prolongation of $A$, which in this introduction we will denote $A^{(n)}$.
Likewise, the $n^{th}$-prolongation
of $X$ is  $P_n(X) = \Spec\  A^{(n)}$.  Thus, there is a bijection between 
differential polynomial functions on $X$ of order $\leq n$ and regular functions on 
$P_n(X)$.  
For each $n$, there is a natural `projection function' $\pi_n:P_n(X) \ra X$, as well as a 
differential polynomial map $\nabla_n:X \ra P_n(X)$, which is a section of $\pi_n$.
In local coordinates, given $\abar \in X$, $\nabla_n(\abar) = (\abar, \dl(\abar), \ldots , \dl^n(\abar))$.
Any differential polynomial function of order $\leq n$ on $X$ factors
as the composition of $\nabla_n$ with a regular function on $P_n(X)$.

One can also define the infinite prolongation $A^{(\infty)}$ 
of $A$ as the ring of all differential polynomial functions on $X$, and the infinite 
prolongation $P_\infty(X)$ as $\Spec\ A^{(\infty)}$.
In this case, $A^{(\infty)}$ is naturally a differential
ring, so that $X^{(\infty)}$ will be what Buium calls a $\DD$-scheme, that is, a scheme over
a differential field, equipped with a sheaf of differential rings.  Note also that
$P_\infty(X)$ is a pro-algebraic variety $P_\infty(X) = \varprojlim P_m(X)$.
As usual, everything above globalizes to arbitrary varieties.

In an appendix to~\cite{Bui93}, Buium notes that prolongations are closely related 
to jet (or arc) spaces of varieties, which have been studied extensively in recent
years (for example, \cite{DL,Cra}).  Recall that the $K$-valued points of the $n^{th}$-jet space
of a variety are the $K[z]/z^{n+1}$-valued points of the variety itself.
Alternatively, the $n^{th}$-jet space represents a certain functor, implicit in
the above characterization.  What Buium observed was that his prolongations 
represent a twisted version of this functor.  Further, given a variety that descends to 
the field of constants, its $n^{th}$-prolongation is isomorphic to its 
$n^{th}$-jet space.  The connection between jets and prolongations does not 
play a large role in Buium's theory, but it is central to the present work.

Below, we develop the theory of prolongations in a rather different way than Buium
does, who built on earlier work of Johnson~\cite{Joh}.
Our approach was inspired by Votja~\cite{Voj}, who gives an elegant construction 
of jet spaces using higher derivations.  The starting point for this paper was the 
observation that, over a differential field, one can modify Vojta's idea
so as to define prolongations in a similar manner.  This perhaps further clarifies the 
relation between jets and prolongations.  It also leads to a rather direct proof
that prolongations represent the twisted jet functor introduced by Buium.  We should also 
mention Gillet's paper~\cite{Gil}, where he develops the theory of prolongations
using adjoint functors, which allows him to give new proofs of 
earlier results of Buium and Kolchin.  

Our work was also motivated by a recent paper of Moosa, Pillay, and Scanlon~\cite{MPS}
on the model theory of differentially closed fields in characteristic 0 with finitely
many commuting derivations.  In that paper, the prolongation of an algebraic variety 
is actually defined in terms of the twisted jet functor.  The authors then go on to 
define more generally prolongations of differential algebraic varieties, which 
are not treated here.  We hope that our paper could be read helpfully as a 
companion to theirs.
Good references for the model theory of differential fields
include \cite{Mar96,Pil02,Sca02}.  For applications to diophantine geometry,
see, for example,~\cite{HP00,PZ,Pil04}.

Let us also say something about higher derivations.  These are defined
below, but the basic idea is that a ring can be equipped with a sequence of additive
maps, $(D_0, D_1, \ldots , D_m, \ldots)$, with $D_1$ a derivation and each 
$D_m$ something like the $m^{th}$ derivative.  In characteristic 0, they are 
essentially equivalent to ordinary derivations, but in positive characteristic, 
higher derivations  are
more general and rather natural to use, for example, for developing
differential Galois theory \cite{MvP}.  From a technical point of view, also,
it was more straightforward to adapt Vojta's construction to fields equipped with 
higher derivations.

\subsection*{Summary of results}
In Section 1, we begin by recalling the definition of a higher derivation (of order $m$)
from an $R$-algebra $A$ to an $R$-algebra $B$.  We then introduce the notion of a
higher derivation over a differential ring $(R,\Dun)$, and show that there is a universal object
$\HS^m_{A/(R,\Dun)}$, which is analogous to the module of K\"{a}hler differentials 
$\Omega_{A/R}$ in the usual case.
This is the $m^{th}$-prolongation of $A$, and we establish some basic properties of it.
(In the first three sections, we restrict our attention to fields with a single higher derivation.
In Section 4, we explain more briefly how to generalize our results to fields with
commuting derivations.)

In Section 2, we define the $m^{th}$-prolongation of a variety, which is by this point
straightforward.  We then prove a characterization of prolongations in terms of 
representable functors,
a result due to Buium in characteristic 0.  We also give a new proof of Moosa, Pillay,
and Scanlon's theorem that the jet space functors and prolongation space functors
commute.

In Section 3, we develop the foundations of Buium's theory of $\DD$-schemes
in arbitrary characteristic.  Here we also prove the main result of this paper,
that the $m^{th}$-prolongation of a variety is isomorphic to its $m^{th}$-jet space.
Previously, this had only been known for $m = 1$, where the $1^{st}$-jet space
is just the tangent variety.  

Finally, in the last section of the paper, we develop the theory of prolongations 
over fields with commuting derivations, very much following the treatment
in the earlier sections of the paper.  

\subsection*{Conventions}
Let $\N = \{0, 1, 2, \ldots \}$ denote the set of natural numbers,
and $\N^+ = \N \setminus \{ 0 \}$.  Throughout the paper, the variable
$m$ will range over the ordered set $\N \cup \{\infty\}$, with $n < \infty$,
for all $n \in \N$.  Variables $i,j,k,l$ will range over $\N$.
We frequently write, for example, $i \leq m$, as shorthand for 
$i \leq m$, if $m \in \N$, and $i < m$, if $m = \infty$.  Likewise,
$i = 0, \ldots, m$ should be taken to mean exactly that, if $m \in \N$,
and to mean $i = 0, 1, \ldots $, if $m = \infty$.  We hope that such shorthand 
will not lead to any unclarity in the presentation.

All rings are commutative with unit.

\section{Higher derivations}
\bdf (See \cite{Mat} or \cite{Voj}.) 
\label{df:HS} 
Let $R$ be a ring, $f: R \ra A$ and 
$R \ra B$ be $R$-algebras, and $m \in \N \cup \{ \infty \}$.
A \emph{higher derivation of order $m$} from $A$ to $B$ over $R$ is a 
sequence $\Dun = (D_0, \ldots , D_m)$, or $(D_0, D_1, \ldots )$ if $m = \infty$,
where $D_0: A \ra B$ is an $R$-algebra homomorphism and $D_1, \ldots , D_m:
A \ra B$ are homomorphisms of (additive) abelian groups such that 
\begin{enumerate}
\item $D_i(f(r))=0$ for all $r\in R$ and $i\geq 1$;
\item (Leibniz Rule) for all $a,b \in A$ and $k \leq m$,
$$D_k(ab)=\sum_{i+j=k}D_i(a)D_j(b).$$
\end{enumerate}
Let $\Der^m_R(A,B)$ denote the set of such derivations.  
\edf  

Higher derivations are also called \emph{Hasse-Schmidt derivations}.

Instead of condition $(1)$, Matsumura requires the $D_i$ to be $R$-module
homomorphisms, which is equivalent.
Below, we will write $(D_0, \ldots, D_m)$, etc.\ even when $m = \infty$.

\bex
With $R,A,B$, and $m$ as above, if $\mathbb{Q} \sub B \sub A$
and $D$ is a usual derivation from 
$A$ to $B$ over $R$, then $D_i = \frac{1}{i!}D^i, i \leq m$, is a higher 
derivation.
\eex

This is a relative notion of higher derivation.  Viewing rings $A$ and $B$ as
$\Z$-algebras, one also gets an absolute notion.  Write $\Der^m(A,B)$ for 
$\Der^m_\Z(A,B)$.

\brm
\label{useful}
Let $R,A,B$, and $m$ be as above, and $\Dun=(D_0,\ldots,D_m)$ be a sequence
of maps from $A$ to $B$.  There is an equivalent condition for $\Dun$ to
be a higher derivation that will be useful below.

For $m < \infty$, let $B_m = B[t]/(t^{m+1})$, the truncated polynomial ring,
and for $m=\infty$, let $B_m= B[[t]]$, the ring of power series.  
It is easy to check that $\Dun$ is a higher derivation if and only if
the map $g:A \ra B_m$,
$$
a \mapsto D_0(a) + D_1(a)t + \ldots + D_m(a)t^m
$$
is a homomorphism of $R$-algebras.

\msk

In the special case $R = \C$, $A=\C(z)$, $B=\C$, and $D_i = 
\frac{1}{i!}d^i/dz$,
this says that the map taking a function $f(z)\in A$ to its 
(truncated) Taylor series expansion around 0 is a $\C$-algebra homomorphism.  
\erm

Observe that $\Der^m_R(A, -)$ is a covariant functor from ($R$-algebras)
to (Sets), and is represented by a (graded)  
$R$-algebra that Vojta calls $\HS^m_{A/R}$,
which is also an $A$-algebra.  (See also Remark~\ref{basics}.(1) below.)
For $m=1$, $\HS^m_{A/R}$ is just the symmetric algebra on $\Om_{A/R}$.

\bdf
Let $A$ be an $R$-algebra.  A \emph{higher derivation} on $A$ is a sequence of 
maps $\Dun \in \Der^\infty_R(A,A)$ such that $D_0 = \textup{Id}_A$.
In this case, we call $(A,\Dun)$ a
{\em \dring\ over} $R$.  
A \emph{homomorphism} $f:(A,\Dun) \ra (B, \Dun)$ between $\DD$-rings over $R$ is an 
$R$-algebra homomorphism such that $f(D_i(a))=D_i(f(a))$, for all $a \in A$ and all $i$.
We will often be interested in the case when $A$
is just a ring (that is, a $\Z$-algebra) and call $(A,\Dun)$ simply a \dring.
The set of \emph{constants} of $(A,\Dun)$ are those $a\in A$ such that
$D_ia=0$, for all $i \geq 1$.
Given $(A,\Dun)$, say that $\Dun$ is \emph{iterative} if for all $i,j$,
$D_i \circ D_j= \tbinom{i+j}{j}D_{i+j}$.
\edf
  
Below we will only consider iterative \drings.  Note that it does not make
sense to talk of a higher derivation from one ring to another being iterative.
  
\brm
As above, let $A_m = A[t]/(t^{m+1})$, for $m<\infty$, and $A_\infty = A[[t]]$.
Let $h:A_m\ra A$ be the homomorphism sending $f(t)\in A_m$ to $f(0)\in A$.
Then a sequence of maps $\Dun = (D_0, \ldots, D_m)$ from $A$ to $A$ is a
higher derivation if and only if the map $d:A \ra A_m$, with $d(a) = \sum_{i\leq m} 
D_i(a)t^i$, is a homomorphism and $h\circ d = \textup{Id}_A$.
\erm

Derivations on a ring extend uniquely to localizations
(see~\cite{Mat}, Theorem 27.2, or \cite{Oku}, Section 1.6, Theorem 1).

\bla
[Quotient Rule]
Let $(A,\Dun)$ be a \dring.  For all invertible $b \in A$, and
$n\in \N^+$,
$$
D_n\left(\frac{1}{b}\right)=\frac{-1}{b}
\left(\sum_{i<n}D_i(b)\cdot D_{n-i}\left(\frac{1}{b}\right)\right).
$$
\ela
\noindent
To obtain this, observe $0=D_n(1)=D_n(b\cdot\frac{1}{b})=
\sum_{i\leq n}D_i(b)\cdot D_{n-i}(\frac{1}{b}))$,
and solve for $D_n(\frac{1}{b})$.

\bla
\label{derloc}
Let $(R,\Dun)$ be a \dring, and  $S$ a multiplicative
subset of $R$.  Then there is a unique extension of $\Dun$ to $S^{-1}R$.
\ela

We now introduce higher derivations on $R$-algebras when $(R,\Dun)$ 
is also a \dring.  This is closely related to Buium's prolongations,
where $(R,\delta)$ is a differential ring, $A$ is an $R$-algebra, and 
one considers derivations on $A$ that are `compatible' with $\delta$.
  
\bdf
Let $(R,\Dun)$ be a \dring.  An $R$-algebra $A$, given by $f:R\ra A$,
is an \emph{$(R,\Dun)$-algebra}
if for all $r\in R$ and all $i$, $f(r)=0$ implies $f(D_i(r))$=0;
in other words, $\Ker(f)$ is a $\DD$-ideal.

Let $f:R\ra A$ and $B$ be $(R,\Dun)$-algebras.  
A \emph{higher derivation} from $A$ to $B$ of order $m$ over $(R,\Dun)$ is a sequence 
$\underline{\delta}= (\delta_0 ,\ldots ,\delta_m)$ such that $\delta_0:A \ra B$
is an $R$-algebra homomorphism, 
$\delta_i: A \ra B, 1 \leq i \leq m$, are (additive) abelian group homomorphisms,
and
\begin{enumerate}
\item $\delta_i(f(r))=\delta_0(f(D_i(r)))$, for $r \in R$;
\item (Leibniz Rule) $\delta_k(ab)=\sum_{i+j=k}\delta_i(a)\delta_j(b)$,
for $a,b\in A$. 
\end{enumerate}
Let $\Der^m_{(R,\Dun)}(A,B)$ denote the set of such derivations.
\edf

Note that if $(R,\Dun)$ is trivial, that is, $D_0 = \textup{Id}_R$ 
and $D_i=0$, $i \geq 1$, then this reduces to Definition~\ref{df:HS}.

As above, given an $(R,\Dun)$-algebra $A$, $\Der^m_{(R,\Dun)}(A, -)$ is a
covariant functor from ($(R,\Dun)$-algebras) to (Sets), which we will now 
observe to be representable.

\bdf
\label{HS}
Let $(R,\Dun)$ be a \dring, $f:R \ra A$ an $(R,\Dun)$-algebra.
For all $m$, define $\HS^m_{A/(R,\Dun)}$ to be the $A$-algebra that is 
the quotient of the polynomial algebra $A[x^{(i)}]_{x\in A, 1\leq i \leq m}$
by the ideal $I$ generated by:
\begin{enumerate}
\item  $(x+y)^{(i)} - x^{(i)} - y^{(i)}: x,y \in A, i=1,\ldots, m;$
\item  $(xy)^{(k)}-\sum_{i+j=k}x^{(i)}y^{(j)}: x,y\in A, k = 1, \ldots , m;$
\item  $f(r)^{(i)}-f(D_i(r)): r\in R, i=1, \ldots ,m.$
\end{enumerate}
In $A[x^{(i)}]$, we identify $x \in A$ with $x^{(0)}$.
There is a universal derivation $\dun=(d_0,\ldots , d_m):A\ra 
\HS^m_{A/(R,\Dun)}$
such that for $i \leq m$ and $x\in A, d_i(x)=x^{(i)}$. 
\edf

\brms
\label{basics}
\begin{enumerate}
\item  With the above notation, if $(R,\Dun)$ is a trivial \dring, then 
$\HS^m_{A/(R,\Dun)}$ is the same as Vojta's $\HS^m_{A/R}$.  In general, 
though, 
$\HS^m_{A/(R,\Dun)}$ is not naturally graded, because of condition $(3)$.
\item  For $m = 1$, we get the first prolongation in the sense of Buium.
\item  For $0\leq m < n \leq \infty$, there are natural $A$-algebra
homomorphisms $f_{mn}:\HS^m_{A/(R,\Dun)} \ra \HS^n_{A/(R,\Dun)}$.
These form a directed system, and 
$$
\HS^\infty_{A/(R,\Dun)}= \lim_{\stackrel{\lra}{i \in \N}} \HS^i_{A/(R,\Dun)}.
$$
\end{enumerate}
\erms

\bdf
Let $(R,\Dun)$ be a \dring.  A $\DD$-$(R,\Dun)$-\emph{algebra} is a \dring\
$(A,\Dun)$ that is also an $(R,\Dun)$-algebra via some map $f:R\ra A$,
such that the derivation on $A$ is compatible with that on $R$.  
That is, for all $r \in R$  and all $i$, $D_i(f(r))=f(D_i(r))$.
\edf

\bla
Given an $(R,\Dun)$-algebra $A$, there is a canonical way to make
$\HS^\infty_{A/(R,\Dun)}$ into a $\DD$-$(R,\Dun)$-algebra.
\ela

\bprf
Extend the universal derivation $\dun:A \ra \HS^\infty_{A/(R,\Dun)}$ to an 
(iterative) higher derivation 
on $\HS^\infty_{A/(R,\Dun)}$ by setting
$$
d_i\left(x^{(j)}\right) = \tbinom{i+j}{i}x^{(i+j)}.
$$
\eprf

\bdf
\label{stuff}
Let $(R,\Dun)$ be a \dring.  Let 
$$
R_m=
\left\{ 
\barr{ll}
R[t]/(t^{m+1}) & \mbox{for } m < \infty \\
R[[t]] & \mbox{for } m = \infty
\earr
\right.
$$
For each $m$, we define a `twisted' homomorphism $e:R \ra R_m$
by $e(r)=D_0(r)+D_1(r)t + \ldots + D_m(r)t^m$.  
Let $\tilde{R}_m$ be the $R$-algebra isomorphic to $R_m$ as a ring,
and made into an $R$-algebra via the map $e:R \ra R_m$.

Let $f:R \ra B$ be an $(R,\Dun)$-algebra.  Define
$B_m = B[t]/(t^{m+1})$, for $m<\infty$, and $B_\infty = B[[t]]$.
Let $\tilde{B}_m$ be the ring $B_m$ made into an $R$-algebra via the map
$\tilde{f}:R \ra \tilde{B}_m$ that sends 
\[
r \mapsto f(D_0(r))+f(D_1(r))t + \ldots + f(D_m(r))t^m.
\]
\edf
  
\bpro
\label{twist}
Let $(R,\Dun)$ be a \dring.   For all $m$, $R_m$ and $\tilde{R}_m$
are isomorphic as $R$-algebras.
\epro

\bprf
Suppose first that $m<\infty$.  We claim that the map $\psi:R_m \ra\tilde{R}_m$
with $\psi(r) = e(r)=D_0(r)+D_1(r)t + \ldots + D_m(r)t^m$, for $r \in R$, and 
$\psi(t) = t$ is an isomorphism of $R$-algebras.  Clearly,
$\psi$ is a homomorphism, so it suffices to check that it is injective and 
surjective.

Let $a=a_0+ a_1t+\ldots +a_mt^m$, so 
$\psi(a)= e(a_0)+e(a_1)t +\ldots + e(a_m)t^m$.  Rearranging terms, 
one gets 
$$  
\psi(a)  =  a_0 
  +  (D_1(a_0) + a_1)t
  +  (D_2(a_0)+D_1(a_1)+a_2)t^2
  +  \cdots
  +  (D_m(a_0) + \ldots + a_m)t^m.
$$  
Suppose that $\psi(a)=0$, so in particular each coefficient of $\psi(a)$
as a polynomial in $t$ is 0.  Thus, $a_0=0$.  Looking at the next term,
$0=D_1(a_0)+a_1=a_1$.  Continuing this way, one sees that all of the 
$a_i$'s are 0, so $a$ itself is 0 and $\psi$ is injective.

To show that $\psi$ is surjective, it suffices to show that for each 
$r \in R$, $r = r + 0t + \ldots + 0t^m\in \tilde{R}_m$ is in Im$(\psi)$.
(Of course, $r \neq \psi(r)$.)  For fixed $r$, we iteratively 
define a sequence,
$c_0, c_1, \ldots , c_m$, of elements of $R_m$ with the following
properties.  One, for all $i \leq m$, the constant term of $\psi(c_i)$, as a 
polynomial in $t$, is $r$.  Two, for $i \geq 1$,  and $1 \leq j \leq i$,
the coefficient of $t^j$ in $\psi(c_i)$ is 0.  Then $\psi(c_m)=r$, as
desired.  Set $c_0=r$.  For the iterative step, suppose that 
$c_0, \ldots ,c_i$ have been defined, and that $\psi(c_i) = 
r+ a_{i+1}t^{i+1} + \ldots + a_mt^m$.  Let $c_{i+1}=
c_i-a_{i+1}t^{i+1}$.  Clearly, this procedure yields such a sequence.

For $m = \infty$, given the isomorphisms $\psi_i:R_i \ra \tilde{R}_i$,
$i<\infty$,
it suffices to note that $R_\infty$ and $\tilde{R}_\infty$ are the
inverse limits of $\{R_i\}_{i<\infty}$ and $\{\tilde{R}_i\}_{i<\infty}$, 
respectively.
The required isomorphism $\psi_\infty:R_\infty\ra\tilde{R}_\infty$
is again given by sending $r\in R$ to $e(r)$, and sending $t$ to $t$.
\eprf

More generally, we have the following.

\bpro
Let $(R,\Dun)$ be a \dring, and $B$ an $(R,\Dun)$-algebra such that
$\Dun$ extends to a derivation on $B$.  Then $\tilde{B}_m\cong B_m$,
as $R$-algebras.
\epro

\bprf
Choose a derivation $\Dun$ on $B$ extending $(R,\Dun)$, and
then argue as above.
\eprf

For fields of characteristic 0, any (higher) derivation on a field $K$
can be extended to a derivation on any extension field $L \sup K$, so 
one has the following corollary.

\bco
Let $(K,\Dun)$ be a \dfield\ of characteristic 0 and $L \sup K$
an extension field.  Then $\tilde{L}_m\cong L_m$,
as $K$-algebras.
\eco

On the other hand, $L_m$ and $\tilde{L}_m$ are not always isomorphic.

\bpro
Let $(K,\Dun)$ be a \dfield\ of characteristic $p > 0$, 
and let $L$ be a purely inseparable 
algebraic extension of $K$, such that there is an $a\in L, b \in K,$
with $a^p=b$ and $D_1(b)\neq 0$.  Then for all $m \geq 1$, 
$L_m$ and $\tilde{L}_m$ are not isomorphic as $K$-algebras.
\epro

\bprf
We show that there is no $K$-algebra homomorphism from $L$ to $\tilde{L}_m$,
which immediately implies the proposition.  In particular, we argue
that any such homomorphism would give an extension of $D_1$ to 
a derivation on $L$, which is impossible  
(as $0\neq D_1(b)=D_1(a^p)=pa^{p-1}D_1(a)=0$, contradiction).

Suppose that $\phi:L\ra \tilde{L}_m$ is a $K$-algebra homomorphism.
For all $c\in K, \phi(c)=D_0(c)+D_1(c)t+ \ldots + D_m(c)t^m$.  
For $x \in L$, write $\phi(x)=\phi_0(x) + \phi_1(x)t + \ldots +
\phi_m(x)t^m$, with $\phi_i:L\ra L$, for $i=0, \ldots, m$.
We claim that for all $x \in L$, $\phi_0(x)=x$.  Indeed, this is clear
for $x \in K$, as $\phi_0(x)=D_0(x)=x$.  Otherwise, $x^{p^n}=y$,
for some $n$ and some $y\in K$.  Then
$$
\phi(x)^{p^n}=\left(\phi_0(x) + \phi_1(x)t + \ldots +\phi_m(x)t^m\right)^{p^n}
=\phi_0(x)^{p^n} +t\cdot g(t)
$$
and also 
$$
\phi(x)^{p^n} = \phi\left(x^{p^n}\right)=\phi(y)=D_0(y)+t\cdot h(t)
$$
with $g(t), h(t)$ polynomials in $L[t]$.
Thus $\phi_0(x)^{p^n}=D_0(y)=y$, so $\phi_0(x)=x$, as desired.

By the claim, $(L,\underline{\phi})$ is a higher derivation of order
$m$ that extends $(K,\Dun)$.  In particular, $\phi_1$ is an extension of 
$D_1$ to $L$, which is impossible.
\eprf

\bpro
\label{mapping}
Let $(R,\Dun)$ be a \dring, $R\ra A$ and $R\ra B$ be $(R,\Dun)$-algebras.
Given a higher derivation $\underline{\dl}=
(\dl_0, \ldots , \dl_m): A \ra B$ there exists a unique 
$(R,\Dun)$-algebra homomorphism, $\phi: \HS^m_{A/(R,\Dun)}\ra B$ such that
$(\dl_0, \ldots, \dl_m)=(\phi\circ d_0, \ldots , \phi \circ d_m).$
Thus $\HS^m_{A/(R,\Dun)}$ (together with the universal derivation 
$\dun : A \ra \HS^m_{A/(R,\Dun)}$) represents the functor $\Der^m_{(R,\Dun)}
(A,-)$.
\epro

\bprf
Define $\phi_0: A[x^{(i)}]_{x\in A, i=1, \ldots , m}\ra B$ by $x^{(i)}
\mapsto \dl_i(x)$.  By the construction of the ideal $I \sub A[x^{(i)}]$
and properties of derivation, we get that $\Ker(\phi_0) \sup I$, so there
is an induced map $\phi:\HS^m_{A/(R,\Dun)}\ra R$.  As $\underline{\dl}
= \phi \circ \dun$, $\phi$ is unique.  Thus the map 
$$
\Hom_R(\HS^m_{A/(R,\Dun)}, B) \lra \Der^m_{(R,\Dun)}(A,B) 
$$ 
is bijective.
\eprf  

Compare the following proposition to Remark~\ref{useful}. 

\bpro
\label{repr}
Let $(R,\Dun)$ be a \dring, $f:R\ra A, g:R\ra B$ be $(R,\Dun)$-algebras.
Given a derivation $\underline{\dl}=(\dl_0, \ldots , \dl_m) \in \Der^m_{(R,\Dun)}(A,B)$,
define a map $\phi = \phi_{\underline{\dl}}:A \ra \tilde{B}_m$ by $\phi(a)=\dl_0(a)+\dl_1(a)t+ \ldots
+\dl_m(x)t^m$.  Then $\phi_{\underline{\dl}} \in \Hom_R(A,\tilde{B}_m)$ and the map
\[
\underline{\dl} \mapsto \phi_{\underline{\dl}}:
\Der^m_{(R,\Dun)}(A,B) \ra \Hom_R(A,\tilde{B}_m)
\]
is a bijection.
\epro

\bprf
The $\dl_i$ are homomorphism of the additive groups, so $\phi$ is also.
The Leibniz Rule implies that $\phi$ is multiplicative, so it only remains to 
show that $\phi\circ f = \tilde{g}$, where $\tilde{g}:R\ra\tilde{B}_m$ is
the homomorphism that makes $\tilde{B}_m$ into an $R$-algebra.  Check,
$$
\barr{lll}
\phi \circ f(x) & = & \dl_0(f(x))+\dl_1(f(x))t+ \ldots + \dl_m(f(x))t^m \\
\\
& = & \dl_0(f(x))+\dl_0(f(D_1(x)))t+\ldots + \dl_0(f(D_m(x)))t^m \\
\\
& = & g(x)+g(D_1(x))t+\ldots+g(D_m(x))t^m\\
\\
& = & \tilde{g}(x)
\earr
$$

This establishes injectivity.  To show surjectivity, we just reverse the 
direction of the argument.  Suppose that $h:A\ra \tilde{B}_m$ is an
$R$-algebra homomorphism, which we can write as
$h(a)=h_0(a) + h_1(a)t + \ldots + h_m(a)t^m$, each $h_i$ a map from 
$A$ to $B$.  We claim that 
$\{h_i: i \leq m\}$ is a higher derivation from $A$ to $B$.  Clearly the 
$h_i$ are additive and satisfy the Leibniz Rule.  So it suffices to show
that for $r\in R$ and $i\leq m$, $h_i(f(r))= h_0(f(D_i(r)))$.  Since $h$
is an $R$-algebra homomorphism, one has $h_i(f(r))= h_0(f(D_i(r)) =g(D_i(r))$.
\eprf

The next corollary follows immediately from Proposition~\ref{mapping}
and Proposition~\ref{repr}.  It is the main point in the characterization of 
prolongations in terms of representable functors.

\bco
[Buium] 
\label{jetdes}
There is a natural bijection
$$
\Hom_R(\HS^m_{A/(R,\Dun)},B) \lra \Hom_R(A,\tilde{B}_m)
$$
\eco  
    
The next result is due to Buium~\cite{Bui93} and Gillet~\cite{Gil}
in a slightly different context.  In fact, Gillet defines the prolongation
functor to be the left adjoint of the forgetful functor
from differential algebras to algebras.

\bpro
\label{adjoint}
Let $(R,\Dun)$ be a \dring, $\Alg_R$ be the category of $(R,\Dun)$-algebras,
and $\DD$-$\Alg_R$ be the category of $\DD$-$(R,\Dun)$-algebras.  
Let $U$ be the forgetful functor
$\DD$-$\Alg_R \ra \Alg_R$.  Then the functor 
$F:\Alg_R \ra \DD$-$\Alg_R$, sending $A$ to $\HS^\infty_{A/(R,\Dun)}$,
is the left adjoint of $U$.
\epro

\bprf
Essentially immediate from the explicit construction given of 
$\HS^\infty_{A/(R,\Dun)}$.   That is, given an $R$-algebra map 
$f: A \ra (B,\Dun^B)$, 
there is an obvious, unique way to lift $f$ to a 
$\DD$-$(R,\Dun)$-algebra map 
$f^\infty: \HS^\infty_{A/(R,\Dun)} \ra (B,\Dun^B)$.  For example, for $x^{(i)}\in\HS^\infty_{A/(R,\Dun)}$,
$x \in A$, then $f^\infty(x^{(i)})=D^B_i(f(x))$.
\eprf

The next result is what Vojta calls the second fundamental exact sequence,
adapted to our context.  For completeness, we include his proof, which
carries over directly.

\bpro
[Second fundamental exact sequence] 
\label{fund}
Let $(R,\Dun)$ be a \dring\ and
$R\ra A \ra B$ a sequence of ring homomorphisms.
Assume that $A\ra B$ is surjective, and let $I$ be its kernel.
Let $J$ be the ideal in $\HS^m_{A/(R,\Dun)}$ generated by
$\{d_ix: i \leq m, x \in I\}$.  Then the following sequence is exact.
\[
0\lra J \lra \HS^m_{A/(R,\Dun)} \lra \HS^m_{B/(R,\Dun)} \lra 0
\]
In the definition of $J$, it suffices to let $x$ vary over a set of 
generators of $I$.
\epro

\bprf
Exactness on the left is immediate.  The natural
map $h:\HS^m_{A/(R,\Dun)} \lra \HS^m_{B/(R,\Dun)}$ is surjective
and its kernel contains $J$, so it remains to show that $\Ker(h)=J$.

From the definition of $\HS^m$, we have the following commutative 
diagram.
$$
\xymatrix{0 \ar[r] & K \ar[r] \ar[d]^{f} &
A[x^{(i)}]_{x\in A , i = 1, \ldots , m} 
\ar[r]\ar[d]^{g}  & \HS^m_{A/(R,\Dun)} \ar[r] \ar[d]^{h} & 0\\
0 \ar[r] & K' \ar[r]  & 
B[x^{(i)}]_{x\in B , i = 1, \ldots , m} \ar[r] &
\HS^m_{B/(R,\Dun)} \ar[r] & 0
}
$$
By Definition~\ref{HS}, the map $f$ is surjective, so by the Snake
Lemma, $\Ker(g)$ maps onto $\Ker(h)$.  But $\Ker(g)$ is 
generated by 
$$
I \cup \{d_ix-d_iy:i=1,\ldots ,m; x,y \in A; x-y \in I\}.
$$
This implies that the kernel of $h$ is generated by
the set $\{d_ix:i=0,\ldots,m,x\in I \}$, as desired.
\eprf

The next two results also occur in Vojta~\cite{Voj}.

\bpro
\label{poly}
Let $(R,\Dun)$ be a \dring, and $A=R[x_i]_{i\in I}$.
Then $\HS^m_{A/(R,\Dun)}$ is the polynomial algebra 
$A[d_jx_i]_{i\in I, j = 1,\ldots, m}$.
\epro

\bprf
Essentially obvious, but also proved in~\cite{Voj}. 
\eprf
 
\bco
\label{polly}
Let $A$ be an $(R,\Dun)$-algebra, $A\cong R[x_i]_{i\in I} / (f_j)_{j\in J}$.
Then 
$$
\HS^m_{A/(R,\Dun)}\cong A[d_kx_i]_{i\in I,k=1,\ldots ,m} / 
(d_kf_j)_{j\in J,k=1,\ldots ,m}.
$$

Suppose further that all of the coefficients of the polynomials 
$f_j,j\in J$, are constants in the ring $R$.  Then
$\HS^m_{A/(R,\Dun)}$ is the same as $\HS^m_{A/R}$, as defined by Vojta.
\eco

\bprf
The first statement follows from Propositions~\ref{fund} and \ref{poly}.  
The second follows from the first, and the analogous statement from~\cite{Voj}. 
\eprf

\section{Prolongations}
\label{oneprol}
In this section, we assume throughout that $(K,\Dun)$ is a \dfield.
Probably everything also works over \drings.
We define prolongations of schemes/varieties over $(K,\Dun)$.  
In characteristic 0, we essentially 
get Buium's prolongations, though there is a slight difference since 
we are using higher derivations.  The construction of the prolongations
is direct, but the results of the previous section provide the connection
with representable functors.  In characteristic 0, this agrees with Buium,
but in characteristic $p>0$, it avoids problems that would arise if one
tries to adapt Buium directly to characteristic $p$, involving 
`dividing by $p$.'
  
\bla
\label{loclz}
Let $A$ be a $(K,\Dun)$-algebra and
and $S$ a multiplicative subset of $A$.
Then there is an isomorphism
$$
\HS^m_{A/(K,\Dun)}\cm_AS^{-1}A \lra \HS^m_{S^{-1}A/(K,\Dun)}.
$$
\ela

\bprf
For all $a\in A$, let $\abar$ denote the image of $a$ in $S^{-1}A$
under the canonical map.  We show that the natural map $\phi$
that sends $d_ia\cm s^{-1}b$ to $s^{-1}b \cdot d_i\abar$,
for all $a,b\in A, s\in S$, and $i\leq m$, is an 
isomorphism.  Clearly, $\phi$ is a homomorphism.
By the quotient rule, for $i \leq m$ and $s^{-1}\bbar\in S^{-1}A$,
$d_i(s^{-1}\bbar)$ can be written as $s^{-n}c$, for some $c \in \HS^m_{A/(K,\Dun)}$,
so $\phi$ is surjective.  To show that $\phi$ is injective,
it suffices to define its inverse.  Let $s^{-1}\bbar \in S^{-1}A$.
We want to define $\phi^{-1}(d_i(s^{-1}\bbar))$ as $c\cm s^{-n}$,
but for this we need to check 
that if $s^{-1}\bbar=t^{-1}\cbar$ in $S^{-1}A$, then 
$\phi^{-1}(d_i(s^{-1}\bbar)) = \phi^{-1}(d_i(t^{-1}\cbar))$,
for all $i \leq m$.
To simplify the presentation, let us assume that $s=t=1$.
As $\bbar = \cbar$ in $S^{-1}A$ is equivalent to there being an $s\in S$
such that $s(c-b)=0$ in $A$, it suffices for us to show that for all 
$a \in A$, if $\abar=0$ in $S^{-1}A$, that is, there is $s\in S$ such that
$sa=0$ in $A$, then $d_ia\cm 1 = 0$ in $\HS^m_{A/(K,\Dun)}\cm_AS^{-1}A$,
for $i \leq m$.

We argue by induction on $i$.  The case $i=0$ is obvious, so assume
that we have proved that $d_ja\cm 1=0$ in $\HS^m_{A/(K,\Dun)}\cm_AS^{-1}A$,
for all $j<i$.  Then
$$
0=sa\cm 1 =d_i(sa)\cm 1 =\sum_{j+k=i}(d_jsd_ka\cm 1)
=sd_ia\cm 1 =(d_ia\cm 1)(1\cm s),
$$
so $(d_ia\cm 1)=0$ in $\HS^m_{A/(K,\Dun)}\cm_AS^{-1}A$, as desired.
\eprf

The next theorem is the differential version of Theorem~4.3 of~\cite{Voj},
and is an easy consequence of Lemma~\ref{loclz}, exactly as in Vojta.

\bthm
\label{prolscheme}
Let $X$ be a $K$-scheme.  For all $m$, there exists a sheaf of 
$\calO_X$-algebras $\HS^m_{X/(K,\Dun)}$ such that (i) for each open affine
$\Spec\,  A\sub X$, there is an isomorphism 
$$
\phi_A:\Gamma(\Spec\, A,\HS^m_{X/(K,\Dun)})\lra\HS^m_{A/(K,\Dun)}
$$
of $(K,\Dun)$-algebras, and (ii) the various $\phi_A$ are compatible with
the localization isomorphism of Lemma~\ref{loclz}.  Moreover, the
collection $((\HS^m_{X/(K,\Dun)}),(\phi_A)_A)$ is unique.
\ethm

\bdf  
Let $X$ be a $K$-scheme.  For all $m$, the 
\emph{$m^{th}$-prolongation of $X$} is the scheme
$$
P_m(X/(K,\Dun)) := \BSpec\, \HS^m_{X/(K,\Dun)}.
$$
Suppose that  $A$ is a $(K,\Dun)$-algebra.  We write $P_m(A/(K,\Dun))=
P_m(\Spec\, A/(K,\Dun))$, which equals $\BSpec\, \HS^m_{A/(K,\Dun)}$.

We will also write $X^{(m)}$ or $P_m(X)$ for $P_m(X/(K,\Dun))$.
\edf
(For the definition of $\BSpec$, see, for example, 
\cite{Har} Ch. II, Ex. 5.17.)

Recall that $K_m = K[t]/(t^{m+1})$, for $m < \infty$,
$K_\infty = K[[t]]$, and that
$e :K \ra K_m$ denotes the twisted homomorphism.
We also let $e:\Spec\, K_m\ra\Spec\, K$ denote the corresponding
twisted morphism of schemes.
Given a $K$-scheme $Y$, let $(Y \times_K\Spec\, K_m)\widetilde{}$ denote the
scheme $(Y \times_K\Spec\, K_m)$ made into a $K$-scheme
via the map $e\circ p:(Y \times_K\Spec\, K_m)\ra \Spec\,  K$,
where $p: (Y \times_K\Spec\, K_m) \ra \Spec\, K_m$ is the canonical projection.
  
\bthm
[Buium]
\label{repthm}
Let $X$ be a $K$-scheme.  For all $m$, the scheme
$P_m(X)$ represents the functor from $K$-schemes to sets given by
$$
Y \mapsto \Hom_K((Y\times_K\Spec\, K_m)\widetilde{}, X).
$$
\ethm

\bprf
For $X$ and $Y$ affine, this follows immediately from Corollary~\ref{jetdes}.
The general case follows by gluing affines.
\eprf

Recall that, given a $K$-scheme $X$, the $m^{th}$-jet space of $X$,
which we denote $J_m(X)$, is the scheme that represents the following 
functor from $K$-schemes to sets.
$$
Y \mapsto \Hom_K(Y\times_K\Spec\, K_m, X)
$$
Buium's theorem clarifies the relationship between prolongations and 
jets.  One also has the following fact, due again to Buium.

\bpro
Let $X$ be a $(K,\Dun)$-scheme such that $X=X'\times_C\Spec\, K$,
where $C$ is the field of constants of $K$, and $X'$ is some $C$-scheme.
(That is, $X$ descends to, or is defined over, $C$.)  Then for all $m$,
$P_m(X)\cong J_m(X)$.
\epro

\bprf
This follows from Corollaries~\ref{polly} and \ref{desc}, below, and
the description of jets in \cite{Voj} 
(see, for example, Theorem 4.3 and Definition 4.4).
\eprf

The next result, due to Moosa, Pillay, and Scanlon, generalizes the well-known 
fact that  for all $m,n \leq\infty$, $J_m(J_n(X))=J_n(J_m(X))$, which 
can be seen by observing that they represent the same functor.
In the original version of~\cite{MPS} it was stated without proof.  A revised
version contains a proof using the Weil restriction.

\bthm
[Moosa, Pillay, and Scanlon]
\label{PJcommute}
Let $X$ be a $K$-scheme.  For all $m,n \leq \infty$,
$$
J_m(P_n(X)) \cong  P_n(J_m(X)).
$$
\ethm
  
\bprf  
We include two proofs.  The first is direct and uses the construction of jets and 
prolongations from \cite{Voj} and this paper.  The second,
closer in spirit  to \cite{MPS}, shows that $J_m(P_n(X))$ and $P_n(J_m(X))$ 
represent the same functor.

It suffices to prove this for affine schemes, so assume that $X=\Spec\, A$.
Even though $K$ is a differential field, we will use $\HS^m_{A/K}$ to denote
the $A$-algebra defined by Vojta, which is defined exactly like 
$\HS^m_{A/(K,\Dun)}$ in Definition~\ref{HS}, 
except that one replaces condition $(3)$ with 
$$
f(r)^{(i)}: r \in K , i = 1, \ldots , m.
$$
The point from our perspective is that $\Spec\, \HS^m_{A/K}$ is the 
$m^{th}$-jet space of $X$, while $\Spec\, \HS^m_{A/(K,\Dun)}$ is the 
$m^{th}$-prolongation of $X$.  Thus
$$
J_m(P_n(X))=\Spec\, \HS^m_{\HS^n_{A/(K,\Dun)}/K}
$$
and
$$
P_n(J_m(X))=\Spec\, \HS^n_{\HS^m_{A/K}/(K,\Dun)}.
$$
We want to show that the two $K$-algebras above are isomorphic.
Let us use $\underline{d}$ for the universal derivation on $\HS^m$,
corresponding to jets, and $\underline{\dl}$ for the universal derivation on $\HS^n$,
corresponding to prolongations.  An arbitrary element of 
$\HS^m_{\HS^n_{A/(K,\Dun)}/K}$ can be written as a sum of terms
$$
d_i\dl_ja: \quad i\leq m,  j\leq n, a \in A,
$$ 
and an arbitrary element of 
$\HS^n_{\HS^m_{A/K}/(K,\Dun)}$ as a sum of terms
$$
\dl_jd_ia:\quad i\leq m, j\leq n, a \in A.
$$  
We claim that the $K$-algebra morphism 
$
\theta: \HS^m_{\HS^n_{A/(K,\Dun)}/K} \ra \HS^n_{\HS^m_{A/K}/(K,\Dun)}
$
with $\theta(d_i\dl_ja) =  \dl_jd_ia$ is an isomorphism.  

First we check that $\theta$ is well-defined.  For example,
$d_i\left(\dl_j(a+b)-\dl_j(a)-\dl_j(b)\right)=0$, so we check 
\[
\theta\left(d_i(\dl_j(a+b)-\dl_j(a)-\dl_j(b))\right)=
\dl_j(d_i(a+b)-d_i(a)-d_i(b))=0.
\]  
Likewise,
$$
\barr{lll}
0 & = & d_i(\dl_j(ab)-\sum_{k+l=j}\dl_k(a)\dl_l(b)\\
& &\\
& = & d_i\dl_j(ab)-\sum_{k+l=j}d_i(\dl_k(a)\dl_l(b))\\
& &\\
& = & d_i\dl_j(ab)-\sum_{k+l=j}\sum_{m+n=i}d_m\dl_k(a)d_n\dl_l(b).
\earr
$$
And we check
$$
\barr{ll}
 &  \theta\left(d_i\dl_j(ab)-\sum_{k+l=j}\sum_{m+n=i}d_m\dl_k(a)d_n\dl_l(b)\right) \\
& \\
=& \dl_jd_i(ab)-\sum_{m+n=i}\sum_{k+l=j}\dl_kd_m(a)\dl_ld_n(b)\\
& \\
=& \dl_jd_i(ab)-\sum_{m+n=i}\dl_j(d_m(a)d_n(b))\\
&\\
=& \dl_j\left(d_i(ab)-\sum_{m+n=i}d_m(a)d_n(b)\right)=0.\\
\earr
$$
Finally, for $c\in K$, we have $0=d_i(\dl_j(c)-\dl_0D_j(c))$.  For $i\geq1$,
$$
\theta(d_i(\dl_j(c)-\dl_0D_j(c)))=\dl_jd_i(c)-\dl_0d_i(D_j(c)) =0
$$
and, for $i=0,$ 
$$
\theta(d_0(\dl_j(c)-\dl_0D_j(c)))=\dl_jd_0(c)-\dl_0d_0(D_j(c))
=\dl_j(c)-\dl_0(D_j(c))=0.
$$

Now that we know that $\theta$ is well-defined, it is clear from the
definition that it respects sums and products.  Finally, it is clearly
a bijection, since there is an obvious inverse.
\msk      

We now give a second proof, along the lines of \cite{MPS}.
Let $Y$ be a $K$-scheme.  There are natural bijections
$$
\barr{ll}
& \Hom_K\left(Y,J_m\left(P_n(X)\right)\right) \\
& \\
\simeq & \Hom_K\left(Y\times_K \Spec\, K_m, P_n(X)\right) \\
& \\
\simeq &\Hom_K
\left(\left(\left(Y\times_K \Spec\, K_m\right)\times_K\Spec\, K_n \right)\widetilde{},X\right)
\earr
$$
and also natural bijections
$$
\barr{ll}
&\Hom_K\left(Y,P_n\left(J_m(X)\right)\right) \\
& \\
\simeq & \Hom_K\left(\left(Y \times_K \Spec\, K_n\right)\widetilde{}, J_m(X)\right) \\
& \\
\simeq & \Hom_K 
\left(\left(\left(Y \times_K \Spec\, K_n\right)\widetilde{} \times_K\Spec\, K_m \right), X \right) \\
\earr
$$
Thus, it suffices to show that for all $Y$,
$$  
\left(\left(Y\times_K \Spec\, K_m\right)\times_K\Spec\, K_n\right)\widetilde{}
\cong
\left(\left(Y \times_K \Spec\, K_n\right)\widetilde{} \times_K\Spec\, K_m\right).
$$

In fact, it suffices to prove this for $Y$ affine.  
We rephrase this as a 
question about isomorphisms of $K$-algebras.  Given a $K$-algebra
$C$, let us write $(C\otimes_K K_n)\widetilde{}$ for what we called
$\tilde{C}_n$ in Definition~\ref{stuff}.  This more closely
parallels our notation for schemes.  That is, for $Z = \Spec\, C$,
then 
$(Z \times_K\Spec\, K_n)\widetilde{} 
= \Spec\left( (C\otimes_K K_n)\widetilde{} \; \right)$.

Everything reduces to showing that, for all 
$Y=\Spec\, B$, the following are isomorphic.
\[
\left((B \otimes_K K_n)\widetilde{} \otimes_KK_m\right) \cong
\left((B\otimes_K K_m)\otimes_KK_n\right)\widetilde{}
\]
Let us write $K_m=K[t]/(t^{m+1})$, $K_n=K[u]/(u^{n+1})$,
and use $e$ for the twisted map from $K$ to $K_n$.

Note that the `trivial'  map
\[
\phi: ((B \otimes_K K_n)\widetilde{}\otimes_KK_m)
\lra
((B\otimes_K K_m)\otimes_KK_n)\widetilde{}
\]
that sends
\[
(1\cm f(t) \cm g(u))
\mapsto
(1\cm g(u)\cm f(t))
\]
is not well-defined.
For example, for $c \in K$, in 
$((B \otimes_K K_n)\widetilde{}\otimes_KK_m)$, 
\[
(1\cm 1 \cm c) = (1\cm e(c) \cm 1)
\]
yet, in $((B\otimes_K K_m)\otimes_KK_n)\widetilde{}$,
\[
\phi(1\cm 1 \cm c) \neq \phi(1\cm e(c) \cm 1).
\]
But a slight variation of this map does work.

First we claim that any non-zero element of $((B \otimes_K K_n)\widetilde{}\otimes_KK_m)$
can be written uniquely as a sum, 
$\sum_{i\leq m, j\leq n}(b_{ij} \cm u^j\cm t^i)$.  Clearly, it suffices to prove
this for elements of the form $(b\cm a_1u^j \cm a_2t^i)$.
And we see that 
\[
\barr{lrll}
&  (b\cm a_1u^j \cm a_2 t^i)  &  =  &  (b\cm e(a_2)a_1u^j \cm t^i)    \\
\\
=  &  \sum_{k\leq n}\left(b\cm D_k(a_2)a_1u^{j+k} \cm t^i \right)  
&  =  &  \sum_{k\leq n} \left(D_k(a_2)a_1b \cm u^{j+k} \cm t^i \right),
\earr
\]
as desired.
Uniqueness is obvious.
Next we observe that this also holds in the algebra 
$((B\otimes_K K_m)\otimes_KK_n)\widetilde{}$.
Note that $(b\cm a_1t^i \cm a_2u^j) \in ((B\otimes_K K_m)\otimes_KK_n)\widetilde{}$
equals $(a_1a_2b\cm t^i\cm u^j)$.

Define 
\[
\theta: ((B \otimes_K K_n)\widetilde{}\otimes_KK_m)
\lra ((B\otimes_K K_m)\otimes_KK_n)\widetilde{}
\]
by
$\theta (b\cm u^j\cm t^i) = (b\cm t^i\cm u^j)$.  
Clearly, $\theta$ is 
a ring homomorphism and injective, but we need to show that it 
is $K$-linear and surjective.  (This sounds completely obvious, but
the $\widetilde{}$\,'s make this more subtle than it first appears.)
Let $c\in K, (b\cm u^j\cm t^i)\in ((B\otimes_K K_n)\widetilde{}\otimes_KK_m)$.
Then 
\[
c\cdot (b\cm u^j\cm t^i)=(b\cm u^j\cm ct^i)=
\sum_{k\leq n}\left(D_k(c)b\cm u^{j+k}\cm t^i\right)
\] 
and 

\[
\barr{llll}
 & \theta \left(\sum_{k\leq n}(D_k(c)b\cm u^{j+k}\cm t^i)\right)
&=&
 \sum_{k\leq n}\left(D_k(c)b \cm t^i \cm u^{j+k}\right) \\
 \\
= &
\sum_{k\leq n}\left(b \cm t^i \cm D_k(c)u^{j+k}\right)          
&=& (b \cm t^i \cm e(c)u^j) \\
\\
= & c\cdot(b \cm t^i \cm u^j).
\earr
\]
This proves $K$-linearity.

To prove that $\theta$ is surjective, it will suffice to show that
for all $c \in K$, that $(1\cm 1 \cm c)\in ((B\otimes_K K_m)\otimes_KK_n)\widetilde{}$
is in the image of $\theta$.  The rest then follows easily.
By Proposition~\ref{twist}, we can rewrite $c$ as
$c =  \sum_{k\leq n}e(c_k)u^k$, so we get that
\[
(1\cm 1 \cm c)= \left(1 \cm 1 \cm \sum_{k\leq n}e(c_k)u^k\right)
= \sum_{k\leq n}\left(e(c_k) \cm 1 \cm u^k\right).
\]  
Thus
$\theta\left(\sum_{k\leq n}(e(c_k)\cm u^k \cm 1)\right) = (1\cm 1 \cm c)$.
\eprf

For completeness, we mention the following, which can be proved 
in the same way as the previous theorem.
\bthm
Let $X$ be a $K$-scheme, $m,n\leq\infty$.  Then
$$
P_m(P_n(X)) = P_n(P_m(X)).
$$
\ethm

\brm
\label{prolsys}
Let $X$ be a $K$-scheme.  For $0\leq m \leq n \leq \infty$,
the maps $f_{mn}:\HS^m_{A/(K,\Dun)} \ra \HS^n_{A/(K,\Dun)}$ of 
Remark~\ref{basics}.(3) give rise to morphisms
$$
f_{mn}:\HS^m_{X/(K,\Dun)}\lra\HS^n_{X/(K,\Dun)}
$$
which again form a directed system.
    
In terms of schemes, the $f_{mn}$ give morphisms
$$
\pi_{nm}:P_n(X/(K,\Dun)) \lra P_m(X/(K,\Dun))
$$
which also form a directed system.  By Remark~\ref{basics}.(3),
$$
\HS^\infty_{X/(K,\Dun)} =
\lim_{\stackrel{\lra}{i\in \N}} \HS^i_{X/(K,\Dun)}
$$
and  
$$
P_\infty(X/(K,\Dun)) =
\lim_{\stackrel{\lla}{i\in \N}} P_i(X/(K,\Dun)).
$$
\erm

\subsection*{Functorial properties}
\label{functor}
There are numerous easy to verify functorial properties of these constructions,
exactly as in \cite{Voj}.  We only mention a few here.
Some more general results hold.

For all $m$,
$\HS^m_{A/(K,\Dun)}$ is functorial in pairs
$(K,\Dun) \ra A$, and $\HS^m_{X/(K,\Dun)}$ and $P_m(X/(K,\Dun))$ are 
functorial in pairs $X \ra \Spec\, K$.  Given a commutative diagram
$$
\xymatrix{
A \ar[r]^\phi & A' \\
(K,\Dun) \ar[u]\ar[r] & (K',\Dun)\ar[u]
}  
$$
there is an induced commutative diagram
$$
\xymatrix{
\HS^m_{A/(K,\Dun)} \ar [r]^{\HS^m_\phi} &  \HS^m_{A'/(K',\Dun)} \\
A \ar[u]^{\phi}\ar[r] & A'\ar[u]
}
$$
that takes $d_ia\in\HS^m_{A/(K,\Dun)}$ to $d_i\phi(a) \in \HS^m_{A'/(K',\Dun)}$
for all $a \in A$ and all $i \leq m$.

Two important cases are base change in $(K,\Dun)$, and functoriality in 
$A$, when $K = K'$.  One also has the following easy lemma.

\bla
Let $A$ be a $(K,\Dun)$-algebra, $(K',\Dun)$ a $\DD$-extension field of $K$,
and $A'=A\cm_KK'$.  Then $\HS^m_{A'/(K',\Dun)}\cong \HS^m_{A/(K,\Dun)}\cm_KK'$
as $A'$-algebras.
\ela

\bprf
Let $\phi$ be the map from $\HS^m_{A'/(K',\Dun)}$ to $\HS^m_{A/(K,\Dun)}\cm_KK'$
that sends $d_k(a\cm c), k\leq m, a\in A,c\in K$, to 
$\sum_{i+j=k}(d_ia\cm 1)(1\cm D_jc))$.  It is clear
that $\phi$ is an isomorphism.
\eprf

These properties carry over to schemes.  
The next result is an easy corollary of the above lemma.

\bco
\label{desc}
Let $(K,\Dun)$ be a \dfield, and let $(K',\Dun)$ be a \dfield\ extension.
Then for all $K$-schemes $X$ and all $m$, 
$$
P_m(X\times_K\Spec\, K')\cong P_m(X)\times_K\Spec\, K'.
$$
\eco

If $f:X\ra X'$ is a morphism of $K$-schemes, one has the 
following commutative diagram, which lifts $f$ to a map between
prolongations.
$$
\xymatrix{
P_m(X) \ar[d] \ar[r]^{P_m(f)} & P_m(X')\ar[d] \\
X \ar[r]^f & X'
}
$$

\bla
Let $X,X'$ be $K$-schemes, and $f:X\ra X'$ a closed immersion.
Then $P_m(f):P_m(X)\ra P_m(X')$ is also a closed immersion.
\ela

\bprf
It is enough to check locally, on affines, where it 
follows from Proposition~\ref{fund}.
\eprf

The following propositions are versions of standard facts about
jet spaces, and can be proved in the same way (for example, see \cite{Bli03}).

\bpro
\label{etmor}
Let $f:X \ra Y$ be an \'etale morphism of schemes over a \dfield\
$(K,\Dun)$.  Then for all $m$,
$$
P_m(X) \cong X \times_Y P_m(Y).
$$
\epro

\bprf
We argue on the corresponding functor of points.  For any 
$K$-scheme $Z$, 
$$
\barr{rll}
\Hom_K(Z,P_m(X)) & \simeq & \Hom_K((Z\times_K \Spec\, K_m)\widetilde{},X) \\
& & \\
\Hom_K(Z,X\times_Y P_m(Y)) & \simeq & 
\Hom_K(Z,X) \times_{\Hom_K(Z,Y)}\Hom_K(Z,P_m(Y)) \\
& & \\
& \simeq & 
\Hom_K(Z,X) \times_{\Hom_K(Z,Y)}\Hom_K((Z\times_K\Spec\, K_m)\widetilde{},Y) 
\earr
$$

Consider the following diagram.

$$
\xymatrix{ X \ar[r] &  Y \\
Z \ar[u]^\phi \ar[r] &
(Z\times_K \Spec\, K_m)\widetilde{}\ar[u]_{\psi}
}
$$
(with a diagonal arrow $\tau$ from $(Z\times_K \Spec\, K_m)\widetilde{}$ to $X$).
A morphism $\tau$ is a $Z$-valued point of $P_m(X)$, and determines
a pair of morphisms 
\[
(\phi,\psi) \in \Hom_K(Z,X) 
\times_{\Hom_K(Z,Y)}\Hom_K((Z\times_K\Spec\, K_m)\widetilde{},Y),
\]
which determines a $Z$-valued point of $X\times_YP_m(Y)$.  This gives
the canonical map from $P_m(X)$ to $X\times_YP_m(Y)$, which does not 
depend on any properties of the morphism $f$.  

In the other direction, a $Z$-valued point of $X\times_Y P_m(Y)$
corresponds to a pair of morphisms $(\phi,\psi)$ making the above
diagram commute.  By formal \'etaleness, there is a unique $\tau$
completing the diagram.  Thus, the map taking $(\phi,\psi)$ to $\tau$
determines the inverse morphism from $X\times_Y P_m(Y)$ to $P_m(X)$,
as desired.
\eprf

\bpro
Let $X$ be a smooth scheme over the \dfield\ $(K,\Dun)$ of dimension 
$n$.  Then for all $m \in \N$, $P_m(X)$ is an $\A^{nm}$-bundle over $X$.
(That is, $X$ can be covered by open sets $U$ such that 
$P_m(U) \cong U \times_K \A^{nm}$.)
\epro

\bprf
By hypothesis, $X \ra \Spec\, K$ is a smooth map, so, by [EGA],
this implies that there is a covering of $X$ by open sets
$U_i$, such that for all $i$, the following diagram commutes
$$
\xymatrix{U_i \ar[d]\ar[r]^{g_i} & \A^n\ar[d]\\
K \ar[r]^{=} & K
}
$$
and $g_i$ is \'etale.  By the previous proposition, 
$P_m(U_i) \cong U_i \times \A^{nm}$, as desired.
\eprf

By the same argument, one also gets the following.

\bco
Let $X$ be an $n$-dimensional  smooth scheme over the \dfield\ $(K,\Dun)$.
Then for all $m$, $P_{m+1}(X)$ is an $\A^{n}$-bundle over $P_m(X)$.
\eco

\section{$\DD$-Schemes}
\label{dschemes}
Many of the definitions and results in this section are from
\cite{Bui93}.

\bdf
Let $(K,\Dun)$ be a \dfield.  A \emph{$\DD$-scheme} over $(K,\Dun)$ is a
$K$-scheme $X$ such that $\calO_X$ is a structure sheaf of 
$\DD$-$(K,\Dun)$-algebras.
A \emph{morphism} of $\DD$-schemes is a morphism of $K$-schemes such that the 
map $\calO_Y \ra f_*\calO_X$ is a map of sheaves of $\DD$-$(K,\Dun)$-algebras.
\edf
  
\bex
Let $X$ be a $K$-scheme.  Then $P_\infty(X)$ is a $\DD$-scheme.
Given a morphism $f:X\ra Y$ of $K$-schemes, the induced map
$f_\infty:P_\infty(X) \ra P_\infty(Y)$ is a morphism of $\DD$-schemes.
\eex

\bpro
Let $(A,\Dun)$ be an $\DD$-$(K,\Dun)$-algebra.  There exists a $\DD$-scheme 
$X=\DSpec(A,\Dun)$ such that, forgetting the $\DD$-structure on $X$,
$X$ is isomorphic to $\Spec\, A$. 
\epro

\bprf
To show that one can add a $\DD$-structure to $\Spec\, A$, it suffices to show
that the localization of a $\DD$-ring is itself a \dring.  
This is the content of Lemma~\ref{derloc}.
\eprf

One also has the following $\DD$-version of a well-known fact from 
algebraic geometry.  (See \cite{Har}, II. Ex. 2.4 and Prop. II.2.3,
or \cite{EH} Thm. I-40.)

\bpro
\label{globsec}
Let $(A,\Dun)$ be a \dring, and $(X,\calO_X)$ a \dsch.  Then there 
is a bijection:
$$
\chi: \Hom_{\DD-\textup{Sch}}(X,\Spec\, A)\lra\Hom_{\DD-\textup{Ring}}
(A,\Gamma(X,\calO_X)).
$$
\epro

\bprf
In the usual case, given a morphism $f:X\ra\Spec\, A$, and the associated 
map $f^\# :\calO_{\Spec\, A} \ra f_* \calO_X$, one gets a homomorphism 
$A\ra\Gamma(X,\calO_X)$ by taking global sections.  This gives a bijection
$$
\chi:\Hom_{\textup{Sch}}(X,\Spec\, A)\lra\Hom_{\textup{Rings}}(A,\Gamma
(X,\calO_X)).
$$
By definition, if $f:X\ra\DSpec\, A$ is a $\DD$-morphism, then the induced 
homomorphism $A \lra \Gamma(X,\calO_X)$ is a homomorphism of \drings, so one 
has an injection:
$$
\chi_\DD:\Hom_{\DD-\textup{Sch}}(X,\DSpec\, A)\lra\Hom_{\DD-\textup{Ring}}
(A,\Gamma(X,\calO_X)).
$$
To verify surjectivity, it suffices to look carefully at the construction of 
$\chi^{-1}$ in [EH].
\eprf

\brm
\label{afsc}
Let $X\sub \A^n$ be an affine $K$-scheme,
$\Gamma(X,\O_X) = K[x_i]_{i=1,\ldots, n}/(f_j)_{j\in J}$.  For all $m$, 
$P_m(X)$ is the closed subscheme of $\A^{nm}=
\Spec(K[d_kx_i]_{i=1,\ldots, n, k=0, \ldots , m})$ with 
$$
\Gamma(P_m(X),\O_{P_m(X)}) = K[d_kx_i]_{i=1,\ldots, n, k=0, \ldots , m}
/(d_kf_j)_{j\in J, k = 0, \ldots , m}.
$$
(This follows from Proposition~\ref{polly}.)
In particular, for every closed point $(a_1,\ldots, a_n)\in X$,
the point $(D_ka_i)_{i=1,\ldots ,n,k=0,\ldots, m}$ is in $P_m(X)$.
The canonical projection from $P_m(X)$ to $X$ maps
a closed point $(a_{ik})_{i=1,\ldots ,n,k=0,\ldots, m}$ to 
its first  $n$ coordinates, $(a_{i0})_{i=1,\ldots,n}$.
\erm

Next, we define $\DD$-polynomial maps between schemes, which we use
to define a section of the canonical map $\pi_m:P_m(X)\ra X$.

\bpro
Let $(K,\Dun)$ be a \dfield.  The prolongation functor, that takes
a $K$-scheme $X$ to the $\DD$-scheme $P_\infty(X)$, is the right adjoint
to the forgetful functor $Y \mapsto Y^!$ from $\DD$-schemes to $K$-schemes.
\epro

\bprf
In \cite{Bui93}, p. 1405. This also follows easily from Proposition~\ref{adjoint}.
\eprf

Recall that given a $K$-scheme $X$, a $K$-rational point of $X$ is a $K$-scheme
homomorphism from $\Spec\, K$ to $X$.  Likewise, if $X$ is a $\DD$-scheme,
we will say that a $K$-rational point of $X$ is a $\DD$-scheme
homomorphism from $\DSpec\, K$ to $X$.  Of course, a $\DD$-morphism
$f:X\ra Y$ naturally induces a map between their $K$-rational points.
The previous proposition immediately implies that there is a natural bijection 
between $K$-rational points of $X$ and of $P_\infty(X)$.

\bdf
Let $X,Y$ be $K$-schemes, and $f:P_\infty(X)\ra P_\infty(Y)$ be a 
$\DD$-morphism.  The natural bijections 
$$
\chi: \Hom_K(\Spec\, K,X)\lra \Hom_{(K,\Dun)}(\DSpec\, K,P_\infty(X))
$$ 
and 
$$
\zeta: \Hom_K(\Spec\, K,Y)\lra \Hom_{(K,\Dun)}(\DSpec\, K,P_\infty(Y))
$$ 
and the induced map 
$$
\hat{f}:\Hom_{(K,\Dun)}(\DSpec\, K,P_\infty(X))\lra\Hom_{(K,\Dun)}(\DSpec\, K,P_\infty(Y))
$$
 determine a 
(set theoretic) map from $K$-rational points of $X$ to those of $Y$,
given by $\zeta^{-1}\circ\hat{f}\circ\chi$.  

A \emph{$\DD$-polynomial map} from $X$ to $Y$ is a map 
on $K$-rational points 
of the form $\zeta^{-1}\circ\hat{f}\circ\chi$,
for some $\DD$-morphism $f:P_\infty(X)\ra P_\infty(Y)$.

Schemes $X$ and $Y$ are \emph{$\DD$-polynomially isomorphic} if 
there are $\DD$-polynomial maps $f:X \ra Y$ and $g:Y\ra X$
such that $g\circ f = \textup{Id}_X$ and $f\circ g = \textup{Id}_Y$.
\edf

\brm
Let $X=\Spec(K[x_i]_{i\leq n}/(f_j)_{j\in J})$, so that 
\[
P_\infty(X) = \Spec(K[d_kx_i]_{i\leq n,k< \infty}/(d_kf_j)_{j\in J,k< \infty}).
\]
The bijection $\chi$ takes $h \in \Hom_K(\Spec\, K,X)$,
which is determined by $b_i=h(x_i), i\leq n,$
to $H \in \Hom_{(K,\Dun)}(\Spec\, K ,P_\infty(X))$ determined by
$D_kb_i = H(d_kx_i), i \leq n$.
\erm

\bpro
\label{dpoly}
Let $X$ be a $K$-scheme, and $m<\infty$.  There exists a $\DD$-polynomial
map $\nabla_m:X\ra P_m(X)$ that is a section of the canonical projection
$p_m:P_m(X)\ra X$.  

Let $f:X\ra Y$ be a morphism of $K$-schemes.   Considering $f$ and
$P_m(f)$ as maps on $K$-rational points, the following diagram commutes.
$$
\xymatrix{
P_m(X) \ar[r]^{P_m(f)} & P_m(Y)\\
X\ar[u]^{\nabla_m} \ar[r]^f & Y \ar[u]^{\nabla_m}
}
$$
\epro

\bprf
By the adjointness of $P_\infty(-)$ and $(-)^!$, there is a natural
bijection
$$
\Hom_K((P_\infty(X))^!,P_m(X))\simeq \Hom_{(K,\Dun)}(P_\infty(X),P_\infty(P_m(X))).
$$ 
Let $f:P_\infty(X)\ra P_\infty(P_m(X))$ be the $\DD$-morphism corresponding
to the canonical projection from $(P_\infty(X))^!$ to $P_m(X)$,
and let $\nabla_m$ be the $\DD$-polynomial map corresponding to $f$.
We show that $\nabla_m$ has the desired properties.

It suffices to check locally, so suppose that $X$ is given as
$\Spec(K[x_i]_{i\leq n}/(f_j)_{j\in J})$.  By Remark~\ref{afsc}, 
$$
\barr{rcl}
P_m(X) & = & \Spec(K[d_kx_i]_{k\leq m, i \leq n}/(d_kf_j)_{k\leq m, j\in J})\\
& & \\
P_\infty(X)& = &\Spec(K[d_kx_i]_{k < \infty, i \leq n}/(d_kf_j)_{k < \infty, j\in J})\\
& & \\
P_\infty(P_m(X)) & = & \Spec(K[d_ld_kx_i]_{i\leq n,k\leq m,l<\infty}/(g_h)_{h\in H})\\
\earr
$$
where $(g_h)_{h\in H}$ is the ideal generated by 
$(d_ld_kf_j)_{j\in J, k \leq m, l < \infty}$.  The $\DD$-morphism  from 
$P_\infty(X)$ to $P_\infty(P_m(X))$, corresponding to the projection 
morphism from $P_\infty(X)$ to $P_m(X)$ is determined by the $\DD$-algebra
homomorphism
\[
K[d_ld_kx_i]_{i\leq n,k\leq m,l<\infty}/(g_h)_{h\in H}
\lra
K[d_kx_i]_{k < \infty, i \leq n}/(d_kf_j)_{k < \infty, j\in J}
\]
that sends $d_ld_kx_i$ to $\tbinom{k+l}{k}d_{k+l}x_i$.
One can then see that this determines the $\DD$-polynomial map from 
$X$ to $P_m(X)$ that takes the closed point $(a_i)_{i\leq n}$ to 
$(D_ka_i)_{i\leq n, k \leq m}$.  By Remark~\ref{afsc}, this is a section 
of $\pi_m$.
  
Next we argue that $P_m(f)\circ\nabla_X=\nabla_Y\circ f$.  Again, it 
suffices to prove this for affine schemes, so assume that 
$X=\Spec\, K[\xbar]/I$ and $Y=\Spec\, K[\ybar]/J$.  
Let $S=K[\xbar]/I$ and $R=K[\ybar]/J$, and let $f$ also denote
the homomorphism from $R$ to $S$ corresponding to $f:X\ra Y$.
A $K$-rational point of $X$ corresponds to a homomorphism $g$ from $S$ to $K$, 
which is determined by the image of $\xbar$, so we think of a $K$-rational point as
a tuple $\abar = g(\xbar)$ of elements of $K$.  Also $P_m(X)=\Spec\, \HS^m_{S/(K,\Dun)}$ 
is affine, and $\HS^m_{S/(K,\Dun)}$ is generated 
by $(d_kx)_{x\in \xbar, k\leq m}$.  We saw above that 
\[
\nabla_X(\abar )=(\abar, D_1(\abar ), \ldots, D_m(\abar )).
\]
More precisely, $\nabla_X(\abar)$ is the $K$-rational point of $P_m(X)$
that corresponds to the map that sends $d_kx\in \HS^m_{S/(K,\Dun)}$ 
to $D_k(g(x))\in K$, for $x\in \xbar$.

Let $f(\abar) = \bbar \in Y$, 
$\bbar =  (g\circ f(y))_{y\in\ybar}$.  
As above,
$\nabla_Y(\bbar)=(\bbar, D_1(\bbar), \ldots, D_m(\bbar))$.
As a map of $K$-algebras, $P_m(f)$ sends
$d_ky$ to $d_kf(y)$, for $y \in \ybar, k \leq m$.
Thus, 
\[
P_m(f)(\abar, D_1(\abar), \ldots, D_m(\abar))= (\bbar, D_1(\bbar ), \ldots, D_m(\bbar ))
\]
 as desired.
\eprf

The following result is new.  It generalizes the well-known fact that 
the first prolongation of a variety is differentially isomorphic to the 
tangent space.  The standard proof is geometric, using the existence,
for any variety $X$, of a differential section $\nabla :X\ra P_1(X)$,
and fact that $P_1(X)$ is a $TX$-torsor.  In contrast, our proof below is 
completely algebraic, though in remarks 
after the proof we try to provide some geometric intuition.
(Recall that, in general for $m > 1$, $J_m(X)$ is not a group scheme over 
$X$, so $P_m(X)$ is not a torsor under $J_m(X)$.  Thus one cannot
generalize the standard proof.)
  
\bthm
Let $X$ be a $K$-scheme.
\begin{enumerate}
\item $P_\infty(P_m(X))$ and $P_\infty(J_m(X))$ are 
isomorphic as $\DD$-schemes.
\item $P_m(X)$ and $J_m(X)$ are $\DD$-polynomially isomorphic.
\end{enumerate}
\ethm
  
\bprf
Part $(2)$ follows immediately from $(1)$, so it suffices to prove $(1)$.
We first establish this for affine schemes.  The general argument 
follows by gluing.  

Let $X = \Spec\, A$.  Note that 
$$
P_\infty(P_m(X))= \Spec \left(\HS^\infty_{\HS^m_{A/(K,\Dun)}/(K,\Dun)} \right)
$$ 
and
$$
P_\infty(J_m(X))= \Spec \left(\HS^\infty_{\HS^m_{A/K}/(K,\Dun)} \right).
$$
Thus we must show that $\HS^\infty_{\HS^m_{A/(K,\Dun)}/(K,\Dun)}$
and $\HS^\infty_{\HS^m_{A/K}/(K,\Dun)}$ are isomorphic
as $\DD$-$(K,\Dun)$-algebras.  Therefore, the theorem follows from
the following proposition.
  
\bpro
Let $A$ be a $(K,\Dun)$-algebra.  Then 
$$
\left(\HS^\infty_{\HS^m_{A/(K,\Dun)}/(K,\Dun)},\Dun \right)
\cong
\left(\HS^\infty_{\HS^m_{A/K}/(K,\Dun)},\Dun \right).
$$
\epro
    
\bprf[Proof of Proposition.]
We first treat the case $A$ a polynomial ring, $A = K[\xbar]$,
where $\xbar$ is a (possibly infinite) tuple.  Write
$$
\barr{lllll}
R & := & \HS^\infty_{\HS^m_{A/(K,\Dun)}/(K,\Dun)} & \cong & K[d_i\dl_jx]
_{0\leq i < \infty,0\leq j\leq m, x \in \xbar} \\
\\
S & := & \HS^\infty_{\HS^m_{A/K}/(K,\Dun)} & \cong & K[d_i\pd_jx]
_{0\leq i < \infty,0\leq j\leq m, x \in \xbar}.
\earr
$$
(Note that $d_i\dl_jx$ and $d_i\pd_jx$ are individual symbols.
One could just have well written instead $x_{ij}$, but the chosen
notation is more suggestive.  Below, we often write
$\dl_jx$, or $\pd_jx$, for $d_0\dl_jx$, or $d_0\pd_jx$,
since we are thinking of $d_0$ as the `identity map.')
Observe that $R$ and $S$ are also $\DD$-rings, letting $D_l(d_i\dl_jx)=
\tbinom{i+l}{l}d_{i+l}\dl_jx$ in $R$, likewise for $S$.
We often write $D_i\pd_j x$, or $D_i\dl_j x$, for 
$d_i\pd_jx$, or $d_i\dl_jx$.

These rings are obviously isomorphic, but we want to 
construct an isomorphism that we can also use in the general
case, $B = K[\xbar]/I$, $I$ any ideal.

Let $\phi: R \ra S$
be the $K$-algebra homomorphism determined by setting
$$
\phi(d_i\dl_jx)=D_i\left(\sum_{k+l=j}d_k\pd_lx\right)
$$ 
for all 
$i,j$, and $x$.  (Of course, $\phi(c)=c$, for $c\in K$.)
Moreover, it is clear that $\phi$ is actually a 
$\DD$-$(K,\Dun)$-algebra homomorphism.

To prove that $\phi$ is an isomorphism, we define a homomorphism
$\psi:S \ra R$ and show that they are inverses of each other.
Let $\psi$ be the homomorphism determined by, for all $i,j$, and $x$,
$$
\psi(d_i\pd_jx)=D_i\left(\sum_{k+l=j}  (-1)^kD_k\dl_lx\right)
$$
As above, one sees easily that $\psi$ is also a 
$\DD$-$(K,\Dun)$-algebra homomorphism.

First, we show that $\psi\circ \phi = \textup{Id}_R$.  It suffices
to calculate this on the generators, $d_i\dl_jx$.
$$
\barr{lcl}
\psi \circ \phi (d_i\dl_jx) & = & \psi \left(D_i\sum_{k+l=j}d_k\pd_lx \right) \\
\\
& = & D_i \left (\sum_{k+l=j}\psi(d_k\pd_lx) \right)\\
\\
& = & D_i \left( \sum_{k+l=j}D_k \left( \sum_{a+b=l}(-1)^aD_a \dl_bx \right) \right) \\
\\
& = & D_i(\dl_jx) + D_i \left( \sum_{b=0}^{j-1} \left( \sum_{a+k=j-b}(-1)^aD_aD_k\dl_b x \right) \right) \\
\\
& = & d_i\dl_jx + D_i \left (\sum_{b=0}^{j-1} \left( \sum_{a+k=j-b}(-1)^a\tbinom{j-b}{a}D_{j-b}\dl_b x \right) \right) \\
\\
& = & d_i\dl_jx \\
\earr
$$  
using the identity $(1-1)^n=\sum_{i=0}^n(-1)^i\tbinom{n}{i}=0$, for $n = j-b$.
Next, we show that $\phi\circ \psi = \textup{Id}_S$, arguing again 
only on generators.
$$
\barr{lcl}
\phi\circ\psi(d_i\pd_jx)& = & 
\phi \left(D_i\sum_{k+l = j } (-1)^kD_k\dl_lx \right)\\
\\
&= &D_i \left(\sum_{k+l=j}(-1)^kD_k\phi(\dl_lx) \right)\\
\\
&= & D_i \left(\sum_{k+l=j}(-1)^kD_k\sum_{a+b=l}D_a\pd_bx \right)\\
\\
& = & D_i (\pd_jx) + \left (\sum_{b=0}^{j-1}\sum_{k+a=j-b}(-1)^kD_kD_a\pd_b x \right)\\
\\
& = & d_i \pd_j x + \left(\sum_{b=0}^{j-1}\sum_{k=0}^{j-b}(-1)^k\tbinom{j-b}{k}\pd_bx \right)\\
\\
& = & d_i \pd_jx\\
\earr
$$

This completes the proof for $A$ a polynomial ring.  We now consider
the general case $B$ a $K$-algebra, $B = K[\xbar]/I$, $I$ an ideal.
By Corollary~\ref{polly} and the analogous result in Vojta, one gets the 
following description of $U: = \HS^\infty_{\HS^m_{B/(K,\Dun)}/(K,\Dun)}$.
(We change the notation slightly, adding  $d_0$ and $\pd_0$ as
`identity functions.')  
$$
U \cong K[d_i\dl_jx]_{0 \leq i, 0 \leq j \leq m, x \in \xbar} / 
(d_i\dl_j f)_{0 \leq i, 0 \leq j \leq m, f \in I}
$$
For a  polynomial $h \in K[\xbar]$, the expression $d_i\dl_jh$  should be considered
shorthand for an element of  the polynomial ring
$K[d_i\dl_jx]_{0 \leq i, 0 \leq j \leq m, x \in \xbar}$, that
can be specified inductively as follows.
$$
\barr{l}
d_i\dl_jc = \tbinom{i+j}{i}D_{i+j}(c) \textup{ for } c \in K; \\
\\
d_i\dl_j(f + g) = d_i\dl_jf + d_i\dl_jg; \\
\\
d_i\dl_j(fg) = \sum_{s+t=i}\left(\sum_{k+l = j}(d_s\dl_kf)(d_t\dl_l g)\right).
\earr
$$
Likewise, we get the following description of 
$V: = \HS^\infty_{\HS^m_{B/K}/(K,\Dun)}$.
$$
V \cong K[d_i\pd_jx]_{0 \leq i, 0 \leq j \leq m, x \in \xbar} / 
(d_i\pd_j f)_{0 \leq i, 0 \leq j \leq m, f \in I}
$$
In $V$, $d_i\dl_jh$ is defined as in $U$, except that for $c\in K$,
$d_i\pd_j c=0$ for $j >0$, and $d_i\pd_j c = D_ic$ for $j = 0$.

Rings $U$ and $V$ are quotients of $R$ and $S$, defined above,
$U = R/J$ and $V=S/L$, where $J$ and $L$ are the ideals
from the definitions of $U$ and $V$.
We claim that the isomorphism $\phi$ from $R$ to $S$ naturally 
induces an isomorphism from $U$ to $V$.  To prove this, it will suffice
to show that $\phi(J) \sub L$ and $\psi(L) \sub J$.  
This follows immediately from the next two claims.

\msk
\noindent
\textbf{Claim 1.}
For each polynomial $h \in A$ and $d_i\dl_jh \in R$, 
$\phi(d_i\dl_jh)=D_i \left(\sum_{k+l=j}d_k\pd_lh \right)$.

\msk
\noindent
\textbf{Claim 2.}
For each polynomial $h \in A$ and  $d_i\pd_jh \in S$,
$\psi (d_i\pd_jh)=D_i \left(\sum_{k+l=j}(-1)^kD_k \dl_lh \right)$.

\msk
\noindent
{\em Proof of Claim 1.}
Since all maps being considered are
additive, it suffices to consider the case $h$ a monomial,
$h = a \ybar$, $a \in K$ and $\ybar = (y_1, \ldots , y_n)$.
We introduce the following multi-index notation.  A multi-index
$\albar$, ($\bebar, \gabar$, etc.) is a sequence of
non-negative integers, $\albar = (\al_1, \ldots , \al_n)$.
We say that the {\em length} of $\albar$ is $n$, and
write $\sum\albar$ for $\sum_i\al_i$.  

Using the (generalized) product rule, we can now give an explicit
definition, in $R$, of 
$$
d_i\dl_ja\ybar :=
D_i\left(\sum_{k+l=j}\dl_ka\dl_l\ybar\right) 
= D_i\left(\sum_{k+l=j}D_ka\cdot\left(\sum_{|\albar|=n,
\sum\albar = l}\dl_{\al_1}y_1\cdot \ldots \cdot \dl_{\al_n}y_n\right)\right)
$$
Likewise, there is an analogous definition for elements of the ring $S$.

$$
d_i\dl_ja\ybar := D_i\left(a\cdot\left(\sum_{|\albar|=n,
\sum\albar = j}\dl_{\al_1}y_1\cdot \ldots \cdot \dl_{\al_n}y_n\right)\right)
$$

We now calculate $\phi(d_i\dl_jh)$ and $D_i(\sum_{k+l=j}d_k\pd_lh)$
to show they are equal.  Thus
$$
\barr{ll}
 & \phi\left(D_i\dl_ja\ybar\right) \\
\\
= & \phi\left(D_i\left(\sum_{k=0}^jD_ka\cdot\left(\sum_{|\albar|=n,
\sum\albar=j-k}\dl_{\al_1}y_1\cdot\ldots\cdot\dl_{\al_n}y_n\right)\right)\right)\\
\\
 =& D_i \left(\sum_{k=0}^jD_ka\cdot\left(\sum_{|\albar|=n,\sum\albar=j-k}
\phi\left(\dl_{\al_1}y_1\cdot\ldots\cdot\dl_{\al_n}y_n\right)\right)\right)\\
\\
=& D_i\left(\sum_{k=0}^jD_ka\cdot\left(\sum_{|\bebar|=|\gabar|=n,\sum\bebar
+\sum\gabar=j-k}
d_{\be_1}\pd_{\ga_1}y_1\cdot\ldots\cdot d_{\be_n}\pd_{\ga_n}y_n\right)\right)
\earr
$$
and 
$$
\barr{ll}
 & D_i\left(\sum_{l+m=j}d_l\pd_ma\ybar\right) \\
\\
=&D_i\left(\sum_{l=0}^jD_l\left(
a\sum_{|\albar|=n,\sum\albar=j-l}\pd_{\al_1}y_1\cdot\ldots\cdot 
\pd_{\al_n}y_n\right)\right)\\
\\
=& D_i\left(\sum_{l=0}^j\sum_{k+m=l}D_ka\cdot\left(D_m\sum_{|\albar|=n,\sum\albar=j-l}
\pd_{\al_1}y_1\cdot\ldots\cdot \pd_{\al_n}y_n\right)\right)\\
\\
=& D_i\left(\sum_{k=0}^jD_ka\cdot\left(\sum_{|\bebar|=|\gabar|=n,\sum\bebar
+\sum\gabar=j-k}
d_{\be_1}\pd_{\ga_1}y_1\cdot\ldots\cdot d_{\be_n}\pd_{\ga_n}y_n\right)\right).
\earr
$$
This completes the proof of Claim 1.

\msk
\noindent
{\em Proof of Claim 2.}
Again we can assume that $h$ is a monomial, $h=a\ybar$, $a \in K$, and
$\ybar=(y_1, \ldots , y_n)$.  
We want to show that $\psi(d_i\pd_ja\ybar)$ equals 
$D_i\sum_{k+l=j}(-1)^kD_k\dl_l(a\ybar)$.  We calculate
$$
\barr{lll}
& \psi\left(d_i\pd_ja\ybar\right)\\
\\
 = & aD_i\psi\left(\pd_j\ybar\right)=aD_i\psi\left(
\sum_{|\albar|=n,\sum\albar=j}\dl_{\al_1}y_1\cdot\ldots\cdot\dl_{\al_n}y_n\right)\\
\\
= & aD_i\left(\sum_{k+l=j}\sum_{|\bebar|=|\gabar|=n, \sum\bebar=k\sum\albar=l}
(-1)^k
d_{\be_1}\dl_{\ga_1}y_1\cdot\ldots\cdot d_{\be_n}\dl_{\ga_n}y_n\right)\\
\\
= &aD_i\left(\sum_{k+l=j}(-1)^kD_k\dl_l\ybar\right)
\earr
$$
and  

$$
\barr{ll}
& D_i\left(\sum_{k+l=j}(-1)^kD_k\dl_l a\ybar\right)\\
\\
= & D_i\left(\sum_{s+t+u+v=j}(-1)^{s+u}
D_s\dl_ta\cdot D_u\dl_v\ybar\right)\\
\\
=& aD_i\left(\sum_{u+v=j}(-1)^uD_u\dl_v\ybar\right)+
D_i\sum_{u+v<j}\left(\sum_{s+t=j-(u+v)}(-1)^{s+u}D_s\dl_ta\cdot D_u\dl_v\ybar\right)\\
\\
=& aD_i\left(\sum_{u+v=j}(-1)^uD_u\dl_v\ybar\right)+
D_i\sum_{u+v<j}\left((-1)^u\sum_{s+t=j-(u+v)}(-1)^s\tbinom{s+t}{s}D_{s+t}a\cdot
D_u\dl_v\ybar\right)\\
\\
=& aD_i\left(\sum_{u+v=j}(-1)^uD_u\dl_v\ybar\right)+
D_i\sum_{u+v<j}\left((-1)^u\sum_{s=0}^{j-(u+v)}(-1)^s\tbinom{j-(u+v)}{s}D_{j-(u+v)}a\cdot
D_u\dl_v\ybar\right)\\
\\
=& aD_i\left(\sum_{u+v=j}(-1)^uD_u\dl_v\ybar\right).
\earr
$$  
This completes the proof of Claim 2, and of the Proposition. 
\eprf
\noindent
This also completes the proof of the Theorem.
\eprf

\brm
Our original definition of $\psi$ was 
$$
\psi(d_i\pd_jx)=D_i\left(\sum_{k+l=j}\sum_{\pi\in P[k]}(-i)^{|\pi|}D_\pi\dl_lx\right)
$$
where $P[k]$ is the set of ordered partitions $\pi$ of $k$, that is,
$\pi=(a_1, \ldots, a_n)\in (\N^+)^n$, with $\sum^n_{p=1}a_p=k$
and $D_\pi=D_{a_1} \circ \cdots \circ D_{a_n}$.  The length of $\pi$ is denoted $|\pi|$.
This is the formula one finds if one inverts $\phi$ `by hand' on examples with 
$i,j$ small.  Later, we observed that the following lemma yields the definition that we
gave above,
$$
\psi(d_i\pd_jx)=D_i\left(\sum_{k+l=j}  (-1)^kD_k\dl_lx\right).
$$
\erm  

\bdf
Given $k\in \N^+$, define a function, `multinomial',
$\mu:P[k]\ra \N^+$, by
\[
\mu(a_1, \ldots , a_n) = 
\tbinom{k}{a_1,\ldots,a_n}:=\frac{k!}{a_1!\cdot\ldots\cdot a_n!}.
\]
\edf
  
One can easily check that $D_\pi = \mu(\pi)D_k$.

\bla
For all $k \in \N^+$, 
\[
\sum_{\pi\in P[k]}(-1)^{|\pi|}\mu(\pi)=(-1)^k.
\]
\ela

\bprf
By induction on $k$.  It will be helpful to stipulate that 
$P[0]:=\{\emptyset\}$, $|\emptyset|=0$, and $\mu(\emptyset)=1$.  
The case $k=1$ is obvious, so assume the lemma holds up to $k-1$.

Given a partition $\pi\in P[j]$, $\pi=(a_1, \ldots, a_m)$,
and $i \in \N^+$,
let $\pi * (i)$ denote the partition $(a_1, \ldots, a_{m}, i) \in P[j+i]$.
Note that $\mu(\pi * (i))=\tbinom{j+i}{i}\mu(\pi).$
$$
\barr{rrll}
& \sum_{\pi\in P[k]}(-1)^{|\pi|} & = &\sum_{i=1}^k\sum_{\pi\in P[k-i]}
(-1)^{|\pi|+1}\mu(\pi * (i))\\
\\
= & -\sum_{i=1}^k\tbinom{k}{i}\sum_{\pi\in P[k-i]}(-1)^{|\pi|}\mu(\pi) &
= & -\sum_{i=1}^k\tbinom{k}{i}(-1)^{k-i} \\
\\
=& (-1)^{k+1}\sum_{i=1}^k\tbinom{k}{i}(-1)^{i}&=&(-1)^k \\
\earr
$$
\eprf

\brm
We now explain the geometric intuition behind the proof of the preceding
proposition.
Recall that by the characterizations of jet spaces and prolongations
via representable functors, we have the following natural bijections,
$$
\barr{rll}
\Hom_K(\Spec\, K, J_m(X)) & \simeq & \Hom_K(\Spec\, K_m, X) \\
\\
\Hom_K(\Spec\, K , P_m(X)) & \simeq & \Hom_K(\Spec\, \tilde{K}_m, X) 
\earr
$$
where, for example, $\Hom_K(\Spec\, K, J_m(X))$ is the set of $K$-rational
points of $J_m(X)$.
By Proposition~\ref{twist}, there is an isomorphism $\Psi$ from 
$K_m$ to $\tilde{K}_m$
so there is a corresponding isomorphism $\Psi'$ from $\Spec\, \tilde{K}_m$
to $\Spec\, K_m$.  Thus, $\Psi'$ induces a natural bijection from
$\Hom_K(\Spec\, K_m, X$) to $\Hom_K(\Spec\, \tilde{K}_m, X)$,
and thus between the $K$-points of $J_m(X)$ and $P_m(X)$.  
When one computes this map in local coordinates,
one gets the morphisms from the proof of the preceding theorem.

We illustrate this for $X= \A^1 = \Spec\, K[x]$.  First, we reformulate
everything in terms of $K$-algebras.  
We have $J_m(X)=\Spec\, K[\pd_0x, \ldots, \pd_mx]$ and
$P_m(X)=\Spec\, K[\dl_0x, \ldots, \dl_mx]$
and the following bijections.
$$
\barr{rll}
\Hom_K(K[\pd_0x, \ldots, \pd_mx], K) & \simeq & \Hom_K(K[x],K_m)\\
\\
\Hom_K(K[\dl_0x, \ldots, \dl_mx], K) & \simeq & \Hom_K(K[x],\tilde{K}_m)\\
\earr
$$
An $m$-jet, $(a_0, \ldots , a_m)\in J_m(X)$, corresponds to the map
\[
f:x\mapsto a_0 + a_1t+ \ldots + a_mt^m \in \Hom_K(K[x],K_m),
\]
which corresponds to the map $F \in \Hom_K(K[\pd_0x, \ldots, \pd_mx], K)$,
where $F(\pd_ix)=a_i$.  Composing $f$ with the isomorphism
$\Psi: K_m \ra \tilde{K}_m$, one gets the map 
\[
g:x\mapsto a_0 + (D_0a_1 + D_1a_0)t+ \ldots + \left(\sum_{j+k=m}D_ja_k\right)t^m 
\in \Hom_K(K[x],\tilde{K}_m),
\]
which corresponds to $G \in \Hom_K(K[\dl_0x, \ldots, \dl_mx], K$,
where $G(\dl_ix)=\sum_{j+k=i}D_ja_k$.

Thus, the bijection above between $K$-points of $J_m(X)$ 
and $P_m(X)$ sends $(a_0, \ldots , a_m)\in J_m(X)$
to $(a_0, D_0a_1+D_1a_0, \ldots,\sum_{j+k=m}D_ja_k )\in P_m(X)$.
This `differential map' from $J_m(X)$ to $P_m(X)$ corresponds to the 
algebraic morphism from $P_\infty(J_m(X))$ to $P_\infty(P_m(X))$
given in the proof of the above theorem by the map $\phi$.
\erm

\section{Multiple derivations}
We now develop the theory of prolongations over a differential field
with finitely many commuting derivations.  In characteristic $p > 0$,
Okugawa~\cite{Oku} developed differential algebra over fields with 
commuting higher derivations.  More recently, differential Galois
theory for such fields has been investigated by Matzat and van der 
Put~\cite{MvP}.
Ziegler~\cite{Zie} has  shown that the model completion of the theory of $n$ 
commuting Hasse-Schmidt derivations is a definitional expansion of the theory
$\textup{SCF}_{p,n}$, the theory of separably closed fields of characteristic
$p$ and degree of imperfection $n$.
Kolchin~\cite{Kol73} considers differential fields with commuting derivations, 
mostly of characteristic 0.  
Moosa, Pillay, and Scanlon~\cite{MPS}
study the model theory of characteristic 0 differential fields
with $n$ commuting derivations.  Since we will consider rings (fields)
of arbitrary characteristic, the results in this section essentially
apply to all of the above contexts.

Of course, it would have been possible to consider multiple derivations
from the beginning.  But it is easier to see the 
theory developed for one derivation first.  The general theory is then
quite similar.

\bdf
For $n \in \N^+$, a ring with $n$ commuting (higher) derivations is
a ring $R$ and a sequence, $\Dun_1, \ldots , \Dun_n$, 
of iterative derivations on $R$, $\Dun_i=(D_{i,0}, D_{i,1}, \ldots )$,
for $i \leq n$, such that for all $i,j,k,l$,
$D_{i,k}\circ D_{j,l}= D_{j,l}\circ D_{i,k}$. 

We also add symbols for `mixed' derivatives.  An $n$-multi-index
$\albar$ is a sequence $(\al_1, \ldots, \al_n)$ of
non-negative integers.  For each $n$-multi-index $\albar$, we 
add an operator $D_{\albar}$, such that 
$D_{\albar} = D_{1,\al_1} \circ \ldots \circ D_{n,\al_n}$.
So $D_{i,j} = D_{\albar}$, where $\albar$ is the multi-index
with a $j$ in the $i^{th}$ place, and 0's everywhere else.
\edf

A ring with $n$ commuting derivations will be written
$(R,\Dun)$ when there no chance of confusion.
These will also be called $\DD$-rings.

Given an $n$-multi-index $\albar$, the \emph{size} of $\albar$, written $|\albar|$,
is the sum $\al_1 + \ldots + \al_n$.  
We will sometimes write $\albar \leq m$ for $|\albar| \leq m$.
There is also a natural partial 
order on $n$-multi-indices, where $\albar \leq \bebar$ 
if and only if for all 
$i \leq n, \al_i \leq \be_i$.  We also writes $\overline{0}$ for the 
multi-index that is a sequence of 0's.  Note that 
$D_{\overline{0}}=\textup{Id}_R$.

Composition of mixed derivatives is completely 
determined by the iteration rule for each derivation, and the fact that
the derivations commute.

\bla
Let $(R,\Dun)$ be a $\DD$-ring, and let $D_\albar,D_\bebar$ mixed derivatives,
$\albar = (\al_1, \ldots, \al_n)$, $\bebar = (\be_1, \ldots , \be_n)$.
$$
D_\albar \circ D_\bebar = \tbinom{\al_1 + \be_1}{\al_1}
\cdots \tbinom{\al_n + \be_n}{\al_n}D_{\albar + \bebar}
$$
\ela

One also has the following generalization of the Leibniz Rule.

\bpro
Let $(R,\Dun$) be a \dring.  Then for any multi-index $\albar$ and 
any $a,b\in R$,
$$
D_\albar(ab)=\sum_{\bebar + \gabar = \albar}D_\bebar(a)D_\gabar(b)
$$
\epro

\bprf
$$
\barr{lcl}
D_\albar(ab) 
& = & D_{1,\al_1}\circ\ldots\circ D_{n,\al_n}(ab) \\
& & \\
& = & D_{1,\al_1} \circ \cdots \circ D_{n-1,\al_{n-1}} 
\left(\sum_{\be_n+\ga_n=\al_n}D_{n,\be_n}(a)D_{n,\ga_n}(b)\right) \\
& & \\
& = & \sum_{\be_n+\ga_n=\al_n}\left(D_{1,\al_1} \circ \cdots \circ D_{n,\al_n} 
\left(D_{n,\be_n}(a)D_{n,\ga_n}(b)\right)\right) \\
& & \\
& \vdots & \\
& & \\
& = &
\sum_{\be_1+\ga_1=\al_1}\left( \ldots \left(\sum_{\be_n+\ga_n=\al_n} 
\left(D_{1,\be_1}\circ \cdots \circ D_{n,\be_n}(a)\cdot D_{1,\ga_1}
\circ\cdots\circ D_{n,\ga_n}(b)\right)\right)\right)\\
\\
& = & \sum_{\bebar + \gabar = \albar}D_\bebar(a)D_\gabar(b)\\
\earr
$$
\eprf
  
\brm
One sees this rule, for example, when taking Taylor series of holomorphic
functions of $n$ variables.  Given functions $f$ and $g$, the Leibniz Rule 
computes the coefficent of $z^\albar$ in the Taylor series of $fg$ from the 
coefficients of the Taylor series of $f$ and of $g$.
\erm

Commuting derivations behave well under localization.

\bla
let $(R,\Dun)$ be a \dring, with $n$ commuting derivations,
and $S \sub R$ a multiplicative subset.  Then the unique extensions of
each of the derivations on $R$ to $S^{-1}R$ also commute.
\ela

\bprf
See~\cite{Oku}, Section 1.6, Corollary 1.
\eprf

\bdf
Let $(R,\Dun)$ be a \dring.  Given
$(R,\Dun)$-algebras $f:R\ra A$ and $B$, a \emph{higher derivation of 
order $m$} from $A$ to $B$ over $(R,\Dun)$ is a set of maps
$\{D_\albar : \albar \leq m\}$ such that $D_{\overline{0}}$ is an
$R$-algebra homomorphism, the $D_{\albar}$ are (additive) abelian group
homomorphisms, and 
\begin{enumerate}
\item  $D_\albar(f(x))=D_{\overline{0}} (f(D_\albar(x)))$;
\item  (Leibniz Rule)  $D_\albar(ab)=\sum_{\bebar + \gabar = \albar}
D_\bebar(a)D_{\gabar}(b)$, for all $a,b, \in A, \albar\leq m$.
\end{enumerate}
\edf

\bdf
Let $(R,\Dun)$ be a \dring,  $f:R \ra A$ an $(R,\Dun)$-algebra.
Define $\HS^m_{A/(R,\Dun)}$ to be the $A$-algebra that is 
the quotient of the polynomial algebra $A[x^{(\albar)}]_{x\in A, 
0 \neq \albar \leq m}$
by the ideal $I$ generated by:
\begin{enumerate}
\item  $(x+y)^{(\albar)} - x^{(\albar)} - y^{(\albar)}: 
  x,y \in A, 0 \neq \albar \leq m;$
\item  $(xy)^{(\albar)}-\sum_{\bebar + \gabar = \albar}
x^{(\bebar)}y^{(\gabar)}: 
x,y\in A, 0 \neq \albar \leq m;$
\item  $f(r)^{(\albar)}-f(D_\albar(r)): r\in R, 0 \neq \albar \leq m .$
\end{enumerate}

In $A[x^{(\albar)}]$, we identify $x \in A$ with $x^{(\overline{0})}$.
There is a universal derivation $\dun=\{d_\albar : \albar\leq m \}:A\ra 
\HS^m_{A/(R,\Dun)}$
such that for $\albar\leq m$ and $x\in A, d_\albar(x)=x^{(\albar)}$. 
\edf

\brm
As before, for $m = \infty$, because $(R,\Dun)$ is an iterative \dring, 
there is a canonical way to make $\HS^\infty_{A/(R,\Dun)}$ into a 
$\DD$-$(R,\Dun)$-algebra.

Extend $\dun:A \ra \HS^\infty_{A/(R,\Dun)}$ to an (iterative) higher 
derivation 
on $\HS^\infty_{A/(R,\Dun)}$ by letting
$$
d_\albar(x^{(\bebar)}) = \tbinom{\al_1+\be_1}{\al_1}\cdots
\tbinom{\al_n+\be_n}{\al_n}x^{(\albar + \bebar)}.
$$
\erm

\bdf
Let $(R,\Dun)$ be a \dring\ with $n$ commuting derivations.
Let
$$
R_m=
\left\{ 
\barr{ll}
R[t_1, \ldots , t_n]/(t_1, \ldots , t_n)^{m+1} 
& \mbox{for } m < \infty \\
R[[t_1,\ldots,t_n]] & \mbox{for } m = \infty
\earr
\right.
$$
For a multi-index $\albar$, write $t^\albar$ as shorthand for
$t_1^{\al_1}\cdots t_n^{\al_n}$.  For each $m$, we define 
the twisted homomorphism $e:R\ra R_m$ by 
$$
e(r)=\sum_{\albar\leq m}D_\albar (r)t^\albar.
$$
Let $\tilde{R}^m$ be the $R$-algebra isomorphic to $R_m$ as a ring, and made 
into an $R$-algebra via the map $e:R\ra R_m$.

Likewise, given an $(R,\Dun)$-algebra $f:R\ra B$, let
$B_m = B[t_1, \ldots , t_n]/(t_1, \ldots , t_n)^{m+1}$, for $m < \infty$,
and $B_\infty = B[[t_1, \ldots , t_n]]$.  
Define $\tilde{f}:R\ra B_m$ by
$$
\tilde{f}(r)=\sum_{\albar\leq m}f(D_\albar(r))t^\albar.
$$
and let $\tilde{B}_m$ be the $(R,\Dun)$-algebra that is the ring $B_m$
with the map $\tilde{f}:R\ra B_m$.
\edf

That $e$ and $\tilde{f}$ are actually homomorphisms follows immediately 
from the Leibniz Rule.
One also has the following converse, whose proof is immediate.

\bla
Let $R$ be a ring and let $f:R \ra R_m$ be a ring
homomorphism, which we write 
$$
f(b)= \sum_{\albar\leq m}f_\albar(b)t^\albar.
$$
Suppose that $f_{\overline{0}}=\textup{Id}_R$.  Then the maps
$\{f_\albar: \albar \leq m\}$ are a higher derivation on $R$.
\ela

\bpro
\label{multtwist}
Let $(R,\Dun)$ be a \dring\ with $n$ commuting derivations.  For all $m$,
 $R_m$ and $\tilde{R}_m$ are isomorphic as $R$-algebras.
\epro

\bprf
The idea of the proof is the same as for Lemma~\ref{twist}.
We first treat the case $m<\infty$.
Let $\psi:R_m\ra\tilde{R}_m$ be the map 
$\psi(r) =e(r)=\sum_{\albar\leq m}D_\albar (r)t^\albar $, for $r \in R$,
and $\psi(t_i) = t_i$, for $i \leq m$.
This is clearly a homomorphism, so it remains to 
check that $\psi$ is injective and surjective.

Linearly order the $n$-multi-indices of size $\leq m$, $\albar_1, \ldots ,
\albar_k$ such that, for all $i,j \leq k$, $|\albar_i|<|\albar_j|$
implies $i < j$.  Let $b \in R_m$ be $b=\sum_{\albar\leq m}b_{\albar}t^\albar$,
and suppose that $\psi(b)=0$.  We will show $b=0$ by showing that 
each $b_{\albar_i} =0$, by induction on $i$.

$$
\barr{lrll}
& \psi(b) & = & \psi\left(\sum_{\albar\leq m}b_{\albar}t^\albar\right) \\
\\
 = & \sum_{\albar\leq m}\psi(b_\albar)t^\albar
& = & \sum_{\albar \leq m} \left(
\sum_{\bebar + \gabar = \albar} D_\bebar(b_\gabar)t^\albar\right)\\
\earr
$$
By assumption, each coefficient 
$\sum_{\bebar + \gabar = \albar} D_\bebar(b_\bebar)$ of $t^\albar$ is $0$.
For the base case, $\albar_1 = \overline{0}$, the constant term,
that is, the coefficient of $t^{\overline{0}}$, is 
$0=D_{\overline{0}}(b_{\overline{0}})=b_{\overline{0}}$.

By induction, suppose that for all $j \leq i$, $b_{\albar_j}=0$.
The $t^{\albar_{i+1}}$ coefficient of $\psi(b)$ is
$0=\sum_{\bebar + \gabar = \albar_{i+1}} D_\bebar(b_\gabar)
=D_{\overline{0}}(b_{\albar_{i+1}})=b_{\albar_{i+1}}$,
because for $\bebar+\gabar=\albar_{i+1}$, if $\bebar\neq 0$,
then $|\gabar| < \albar_{i+1}$, so $b_\gabar = 0$, by the induction
hypothesis.

To show that $\psi$ is surjective, it suffices to show that for each 
$r \in R$, $r \in \tilde{R}_m$ is in Im$(\psi)$.
For fixed $r$, we iteratively  define a sequence,
$c_0, c_1, \ldots , c_k$, of elements of $R_m$ with the following
properties.  
One, for all $i \leq k$, the constant term of $\psi(c_i)$, as a 
polynomial in the $t_i$, is $r$.  Two, for $i \geq 1$,  and $1 \leq j \leq i$,
the coefficient of $t^{\albar_j}$ in $\psi(c_i)$ is 0.  Then $\psi(c_m)=r$, as
desired.  Set $c_0=r$.  For the iterative step, suppose that 
$c_0, \ldots ,c_i$ have been defined, and that $\psi(c_i) = 
r+ \sum_{i+l \leq j \leq k}a_{\albar_j}t^{\albar_{j}}$.  Let $c_{i+1}=
c_i-a_{\albar_{i+1}}t^{\albar_{i+1}}$.  
Clearly, this procedure yields such a sequence.

For $m = \infty$, given the isomorphisms $\psi_i:R_i \ra \tilde{R}_i$,
$i < \infty$,
it suffices to note again that $R_\infty$ and $\tilde{R}_\infty$ are the
inverse limits of $\{R_i\}_{i<\infty}$ and $\{\tilde{R}_i\}_{i<\infty}$, 
respectively.
The required isomorphism $\psi_\infty:R_\infty\ra\tilde{R}_\infty$
also sends $r\in R$ to $e(r)$, and sends each $t_i$ to $t_i$.
\eprf

The next two results are proved in the same way as Proposition~\ref{mapping}
and Lemma~\ref{repr}, respectively.

\bpro
\label{mapping2}
Let $(R,\Dun)$ be a \dring, and $R\ra A$ and $R\ra B$ be $(R,\Dun)$-algebras.
Given a higher derivation $\underline{\dl}
\in\Der^m_{(R,\Dun)}(A,B)$, there exists a unique 
$(R,\Dun)$-algebra homomorphism, $\phi: \HS^m_{A/(R,\Dun)}\ra B$ such that
for all $\albar \leq m$, $\dl_\albar = \phi \circ d_\albar$.
Thus $\HS^m_{A/(R,\Dun)}$ (together with the universal derivation 
$\dun : A \ra \HS^m_{A/(R,\Dun)}$) represents the functor $\Der^m_{(R,\Dun)}
(A,-)$.
\epro

\bla
\label{repr2}
Let $(R,\Dun), R\ra A, R\ra B$, and $m$ be as above.
Given $\underline{\dl}\in \Der^m_{(R,\Dun)}(A,B)$,
define a map
 $\phi =  \phi_{\underline{\dl}}:A \ra \tilde{B}_m$ 
by 
$\phi(a)=\sum_{\albar\leq m}\dl_\albar(a)t^\albar$.
Then $\phi \in \Hom_R(A,\tilde{B}_m)$ and the map
\[
\underline{\dl} \mapsto \phi_{\underline{\dl}}:
\Der^m_{(R,\Dun)}(A,B) \lra \Hom_R(A,\tilde{B}_m)
\]
is a bijection.
\ela

The next corollary is the key result in characterizing 
prolongations in terms of representable functors, as in Buium.

\bco
There is a natural bijection
$$
\Hom_R(\HS^m_{A/(R,\Dun)},B) \lra \Hom_R(A,\tilde{B}_m).
$$
\eco  

\bprf
Immediate from Proposition~\ref{mapping2} and Lemma~\ref{repr2}.
\eprf

The following results are proved as in the case of a single derivation.

\bpro
\label{multadjoint}
Let $(R,\Dun)$ be a \dring, $\Alg_R$ be the category of $(R,\Dun)$-algebras,
and $\DD$-$\Alg_R$ be the category of $\DD$-$(R,\Dun)$-algebras.
Let $U$ be the forgetful functor $\DD$-$\Alg_R \ra \Alg_R$.  Then the functor 
$F:\Alg_R \ra \DD$-$\Alg_R$, sending $A$ to $\HS^\infty_{A/(R,\Dun)}$,
is the left adjoint of $U$.
\epro

\bpro
[Second fundamental exact sequence] 
\label{multfund}
Let $(R,\Dun)$ be a \dring\ and
$R\ra A \ra B$ a sequence of ring homomorphisms.
Assume that $A\ra B$ is surjective, and let $I$ be its kernel.
Let $J$ be the ideal in $\HS^m_{A/(R,\Dun)}$ generated by
$\{d_\albar x: \albar  \leq m, x \in I\}$.  Then the following sequence is exact.
\[
0\lra J \lra \HS^m_{A/(R,\Dun)} \lra \HS^m_{B/(R,\Dun)} \lra 0
\]
In the definition of $J$, it suffices to let $x$ vary over a set of 
generators of $I$.
\epro

\bpro
\label{multpoly}
Let $(R,\Dun)$ be a \dring, and $A=R[x_i]_{i\in I}$.
Then $\HS^m_{A/(R,\Dun)}$ is the polynomial algebra 
$A[d_\albar x_i]_{i\in I, 1 \leq \albar \leq m}$.
\epro

\brm
\label{dimcalc}
For $m,n \geq 1$, define $c_{n,m}$ to be the number of $n$-multi-indices
of size $\leq m$.  Equivalently, $c_{n,m}$ is the number of monomials
in $n$ variables of order $\leq m$ or the number of mixed partial derivatives 
in $n$ variables of total order $\leq m$.

By the previous proposition, given a polynomial ring $A= R[x_1, \ldots, x_q]$
over a ring $R$ with $n$ commuting derivations, then $\HS^m_{A/(R,\Dun)}$
is a polynomial ring in $q\cdot c_{n,m}$ indeterminates.  
\erm

\bco
\label{multpolly}
Let $A$ be an $(R,\Dun)$-algebra, $A\cong R[x_i]_{i\in I} / (f_j)_{j\in J}$.
Then 
$$
\HS^m_{A/(R,\Dun)}\cong A[d_\albar x_i]_{i\in I, 1 \leq \albar \leq m } / 
(d_\albar f_j)_{j\in J, 1 \leq  \albar \leq m }.
$$
\eco

\subsection{Prolongations}
In this section, we generalize the results of Section~\ref{oneprol} to 
fields with many derivations.  Almost everything goes through as before.
Assume throughout that $(K,\Dun)$ is a $\DD$-field with $n$ commuting derivations.

\bla
\label{multloclz}
Let $A$ be a $(K,\Dun)$-algebra
and $S$ a multiplicative subset of $A$.
There is an isomorphism
$$
\HS^m_{A/(K,\Dun)}\cm_AS^{-1}A \lra \HS^m_{S^{-1}A/(K,\Dun)}.
$$
\ela

\bthm
Let $X$ be a $K$-scheme.  For all $m$, there exists a sheaf of 
$\calO_X$-algebras $\HS^m_{X/(K,\Dun)}$ such that (i) for each open affine
$\Spec\, A\sub X$, there is an isomorphism 
$$
\phi_A:\Gamma(\Spec\, A,\HS^m_{X/(K,\Dun)})\lra\HS^m_{A/(K,\Dun)}
$$
of $(K,\Dun)$-algebras, and (ii) the various $\phi_A$ are compatible with
the localization isomorphism of Lemma~\ref{multloclz}.  Moreover, the
collection $((\HS^m_{X/(K,\Dun)}),(\phi_A)_A)$ is unique.
\ethm

\bdf  
Let $X$ be a $K$-scheme.  For all $m$, the 
{\em $m^{th}$-prolongation of $X$} is the scheme
$$
P_m(X/(K,\Dun)) := \BSpec\, \HS^m_{X/(K,\Dun)}.
$$
Suppose that $A$ is a $(K,\Dun)$-algebra. We write $P_m(A/(K,\Dun))=
P_m(\Spec\, A/(K,\Dun))$, which equals $\BSpec\, \HS^m_{A/(K,\Dun)}$.

We will also write $X^m$ or $P_m(X)$ for $P_m(X/(K,\Dun))$.
\edf

Recall that for $m < \infty$, 
$K_m = K[t_1, \ldots , t_n]/ (t_1, \ldots , t^{m+1})$, $K_\infty = K[[t_1, \ldots t_n]]$, 
and that $e :K \ra K_m$ denotes the twisted homomorphism.
We also let $e:\Spec\, K_m\ra\Spec\, K$ denote the corresponding
twisted morphism of schemes.
Given a $K$-scheme $Y$, let $(Y \times_K\Spec\, K_m)\widetilde{}$ denote the
scheme $(Y \times_K\Spec\, K_m)$ made into a $K$-scheme
via the map $e\circ p:(Y \times_K\Spec\, K_m)\ra \Spec\, K$,
where $p: (Y \times_K\Spec\, K_m) \ra \Spec\, K_m$ is the canonical projection.
  
\bthm
Let $X$ be a $K$-scheme.  For all $m$, the scheme
$P_m(X)$ represents the functor from $K$-schemes to sets given by
$$
Y \mapsto \Hom_K((Y\times_K\Spec\, K_m)\widetilde{}, X).
$$
\ethm

\bthm
[Moosa, Pillay, and Scanlon]
Let $X$ be a $K$-scheme.  For all $m,q \leq \infty$,
$$
J_m(P_q(X)) \cong  P_q(J_m(X)).
$$
\ethm
  
\bprf  
Both proofs of Theorem~\ref{PJcommute} generalize easily.  Here we only show how
to adapt the second proof.  Exactly as before, it suffices to show that for any
$K$-algebra $B$, the following are isomorphic.
$$
\left((B \otimes_K K_q)\widetilde{}\otimes_KK_m\right) \cong
\left((B\otimes_K K_m)\otimes_KK_q\right)\widetilde{}
$$
where $K_m=K[t_1, \ldots t_n]/(t_1, \ldots , t_n)^{m+1}$, 
$K_q=K[u_1, \ldots, u_n]/(u_1, \ldots , u_n)^{q+1}$,
and we use $e$ for the twisted map from $K$ to $K_q$.

We claim any non-zero element of $((B\cm_K K_q)\widetilde{}\cm_KK_m)$ can be written uniquely
as a sum 
\[
\sum_{\albar\leq m, \bebar \leq q}(b_{\albar,\bebar}\cm u^\bebar \cm t^\albar).
\]
Again, it suffices to prove this for elements of the form $(b\cm a_1u^\bebar\cm a_2 t^\albar)$.
And 
\[
\barr{lrll}
& (b\cm a_1u^\bebar\cm a_2 t^\albar) 
&=& 
(b \cm e(a_2)a_1u^\bebar\cm t^\albar)    \\
\\
= &\sum_{\gabar\leq q}(b\cm D_\gabar(a_2)a_1u^{\bebar + \gabar} \cm t^\albar)
& = &
\sum_{\gabar \leq q}(D_\gabar(a_2)a_1b\cm u^{\bebar+ \gabar}\cm t^\albar)
\earr
\]
 as desired.
Secondly, observe that this also holds in $((B\cm_KK_m)\cm_KK_q)$, as
$(b\cm a_1t^\albar\cm a_2u^\bebar)\in ((B\cm_KK_m)\cm_KK_q)$ equals
$(a_1a_2b\cm t^\albar \cm u^\bebar)$.
  
Define 
\[
\theta: ((B\cm_K K_q)\widetilde{}\cm_KK_m) \lra
((B\cm_KK_m)\cm_KK_q)\widetilde{}
\] 
by 
$\theta(b\cm u^\bebar \cm t^\albar) = (b\cm t^\albar \cm u^\bebar)$.  
It suffices to show that $\theta$ is $K$-linear and  surjective.
Let $c \in K$, $(b\cm u^\bebar\cm t^\albar)\in((B\cm_KK_q)\widetilde{}\cm_KK_m)$.
Then 
\[
c\cdot(b\cm u^\bebar\cm t^\albar)
= \sum_{\gabar\leq q}(D_\gabar(c)b\cm u^{\bebar + \gabar} \cm t^\albar),
\] 
and
\[
\barr{llll}
& \theta\left(\sum_{\gabar\leq q}(D\gabar(c)b\cm u^{\bebar + \gabar} \cm t^\albar)\right)
&=&
\sum_{\gabar\leq q}(D_\gabar(c)b\cm t^\albar \cm u^{\bebar + \gabar}) \\
\\
=&
\sum_{\gabar\leq q}(b\cm t^\albar \cm D_\gabar(c)u^{\bebar + \gabar})
&=&
(b\cm t^\albar \cm e(c)u^{\bebar}) \\
\\
=&
 c\cdot(b\cm t^\albar \cm u^\bebar). & &
\earr
\]
This proves $K$-linearity.

To prove that $\theta$ is surjective, it will suffice to show that for all $c \in K$,
that $(1 \cm 1 \cm c) \in ((B\cm_K K_m) \cm_K K_q)\widetilde{}$ is in the image of $\theta$.
The rest follows easily.  By Proposition~\ref{multtwist}, we can write $c$ as 
$c = \sum_{\gabar\leq q}e(c_\gabar)u^\gabar$, so we get that 
\[
 (1 \cm 1 \cm c) =
  \left(1 \cm 1 \cm \sum_{\gabar\leq q}e(c_\gabar)u^\gabar\right)
= \sum_{\gabar\leq q}\left(e(c_\gabar)\cm 1 \cm  u^\gabar\right).
\]
Thus 
\[
\theta\left(\sum_{\gabar\leq q}(e(c_\gabar)\cm  u^\gabar \cm 1)\right)
= (1 \cm 1 \cm c).
\]
\eprf

\brm
Let $X$ be a $K$-scheme.  As in Remark~\ref{prolsys}, 
for $0 \leq m \leq n \leq \infty$, the canonical maps $f_{mn}:
\HS^m_{A/(K,\Dun)} \ra \HS^n_{A/(K,\Dun)}$ determine a directed
system of morphisms
$$
f_{mn}:\HS^m_{X/(K,\Dun)}\lra\HS^n_{X/(K,\Dun)}.
$$
    
In terms of schemes, the $f_{mn}$ give morphisms
$$
\pi_{nm}:P_n(X/(K,\Dun)) \lra P_m(X/(K,\Dun))
$$
which also form a directed system.  Exactly as above, we also have
$$
\HS^\infty_{X/(K,\Dun)} =
\lim_{\stackrel{\lra}{i\in \N}} \HS^i_{X/(K,\Dun)}
$$
and  
$$
P_\infty(X/(K,\Dun)) =
\lim_{\stackrel{\lla}{i\in \N}} P_i(X/(K,\Dun)).
$$
\erm
  
\subsection*{Functorial properties}
There are many functorial properties of these constructions, 
precisely as discussed on page~\pageref{functor}.

\bla
Let $A$ be a $(K,\Dun)$-algebra, $(K',\Dun)$ a $\DD$-extension field of $K$,
and $A'=A\cm_KK'$.  Then $\HS^m_{A'/(K',\Dun)}\cong \HS^m_{A/(K,\Dun)}\cm_KK'$
as $A'$-algebras.
\ela

\bprf
Let $\phi$ be the map from $\HS^m_{A'/(K',\Dun)}$ to $\HS^m_{A/(K,\Dun)}\cm_KK'$
that sends $d_\albar(a\cm c), \albar \leq m, a\in A,c\in K$, to 
$\sum_{\bebar + \gabar =k}(d_\bebar a\cm 1)(1\cm D_\gabar c))$.  It is clear
that $\phi$ is an isomorphism.
\eprf

\bco
Let $(K,\Dun)$ be a \dfield, and let $(K',\Dun)$ be a \dfield\ extension.
Then for all $K$-schemes $X$ and all $m$,
$$
P_m(X\times_K\Spec\, K')\cong P_m(X)\times_K\Spec\, K'.
$$
\eco

As above, if $f:X\ra X'$ is a morphism of $K$-schemes, 
then there is an induced map $P_m(f):P_m(X)\ra P_m(X')$ between
their prolongations.

\bla
Let $X,X'$ be $K$-schemes, and $f:X\ra X'$ a closed immersion.
Then $P_m(f):P_m(X)\ra P_m(X')$ is also a closed immersion.
\ela

\bpro
Let $f:X \ra Y$ be an \'etale morphism of schemes over a \dfield\
$(K,\Dun)$.  Then for all $m$,
$$
P_m(X) \cong X \times_Y P_m(Y).
$$
\epro

\brm
Notice by Remark~\ref{dimcalc} that for any $q$ and any $m < \infty$, 
$\textup{dim}(P_m(\A^q)) = q\cdot c_{n,m}$.  
\erm

\bpro
Let $X$ be a smooth scheme over the \dfield\ $(K,\Dun)$ of dimension 
$q$.  Then for all $m$, $P_m(X)$ is an $\A^{q\cdot c_{n,m}}$-bundle over $X$.
(That is, $X$ can be covered by open sets $U$ such that 
$P_m(U) \cong U \times_K \A^{q\cdot c_{n,m}}$.
\epro

\bprf
By hypothesis, $X \ra \Spec\, K$ is a smooth map, so, by [EGA],
this implies that there is a covering of $X$ by open sets
$U_i$, such that for all $i$, the following diagram commutes
$$
\xymatrix{U_i \ar[d]\ar[r]^{g_i} & \A^n \ar[d]\\
K \ar[r]^{=} & K
}
$$
and $g_i$ is \'etale.  By the previous proposition, 
$P_m(U_i) \cong U_i \times \A^{q\cdot c_{n,m}}$, as desired.
\eprf

\bco
Let $X$ be a smooth scheme over the \dfield\ $(K,\Dun)$ of dimension 
$q$.  Then for all $m$, $P_{m+1}(X)$ is an 
$\A^{q(c_{n,m} - c_{n,m-1})}$-bundle over $P_m(X)$.
\eco

\subsection{$\DD$-Schemes}
\label{multdschemes}
We generalize material from Section~\ref{dschemes}, which is straightforward.

\bdf
Let $(K,\Dun)$ be a \dfield.  A {\em $\DD$-scheme} over $(K,\Dun)$ is a
$K$-scheme $X$ such that $\calO_X$ is a structure sheaf of 
$\DD$-$(K,\Dun)$-algebras.
A {\em morphism} of $\DD$-schemes is a morphism of $R$-schemes such that the 
map $\calO_Y \ra f_*\calO_X$ is a map of sheaves of $(R,\Dun)$-algebras.
\edf

\bpro
Let $(A,\Dun)$ be a $\DD$-$(K,\Dun)$-algebra.  There exists a $\DD$-scheme 
$X=\DSpec(A,\Dun)$ such that, forgetting the $\DD$-structure on $X$,
$X$ is isomorphic to $\Spec\, A$. 
\epro

\bpro
Let $(A,\Dun)$ be a \dring, and $(X,\calO_X)$ a \dsch.  Then there 
is a bijection:
$$
\chi: \Hom_{\DD-\textup{Sch}}(X,\Spec\, A)\lra\Hom_{\DD-\textup{Ring}}
(A,\Gamma(X,\calO_X)).
$$
\epro

\brm
\label{multafsc}
Let $X\sub \A^q$ be an affine $K$-scheme,
$\Gamma(X,\O_X) = K[x_i]_{i=1,\ldots, q}/(f_j)_{j\in J}$.  For all $m\leq\infty$, 
$P_m(X)$ is the closed subscheme of $\A^{q\cdot c_{n,m}}=
\Spec(K[D_\albar x_i]_{i=1,\ldots, q, \albar \leq m})$ with 
$$
\Gamma(P_m(X),\O_{P_m(X)}) = K[D_\albar x_i]_{i=1,\ldots, q, \albar \leq m}
/(D_\albar f_j)_{j\in J, \albar \leq m}.
$$
(This follows from Proposition~\ref{multpolly}.)
In particular, for every closed point $(a_1,\ldots, a_q)\in X$,
the point $(D_\albar a_i)_{i=1,\ldots ,q,\albar \leq  m}$ is in $P_m(X)$.
The canonical projection from $P_m(X)$ to $X$ maps
a closed point $(a_{i,\albar})_{i=1,\ldots ,q,\albar \leq  m}$ to 
its first  $q$ coordinates, $(a_{i,0})_{i=1,\ldots,q}$.
\erm

\bpro
Let $(K,\Dun)$ be a \dfield.  The prolongation functor, that takes
a $K$-scheme $X$ to the $\DD$-scheme $P_\infty(X)$, is the right adjoint
to the forgetful functor $Y \mapsto Y^!$ from $\DD$-schemes to $K$-schemes.
\epro

As before, if $X$ is a $\DD$-scheme,
we define a $K$-rational point of $X$ to be a $\DD$-scheme
homomorphism from $\DSpec\, K$ to $X$.  Of course, a $\DD$-morphism
$f:X\ra Y$ naturally induces a map between their $K$-rational points.
The previous proposition immediately implies that there is a natural bijection 
between $K$-rational points of $X$ and of $P_\infty(X)$.

\bdf
Let $X,Y$ be $K$-schemes, and $f:P_\infty(X)\ra P_\infty(Y)$ be a 
$\DD$-morphism.  The natural bijections 
$$
\chi: \Hom_K(\Spec\, K,X)\lra \Hom_{(K,\Dun)}(\DSpec\, K,P_\infty(X))
$$ 
and 
$$
\zeta: \Hom_K(\Spec\, K,Y)\lra \Hom_{(K,\Dun)}(\DSpec\, K,P_\infty(Y))
$$ 
and the induced map 
$$
\hat{f}:\Hom_{(K,\Dun)}(\DSpec\, K,P_\infty(X))\lra\Hom_{(K,\Dun)}(\DSpec\, K,P_\infty(Y))
$$
determine a (set theoretic) map from $K$-rational points of $X$ to those of $Y$,
given by $\zeta^{-1}\circ\hat{f}\circ\chi$.  

A {\em $\DD$-polynomial map} from $X$ to $Y$ is a map 
on $K$-rational points 
of the form $\zeta^{-1}\circ\hat{f}\circ\chi$,
for some $\DD$-morphism $f:P_\infty(X)\ra P_\infty(Y)$.

Schemes $X$ and $Y$ are {\em $\DD$-polynomially isomorphic} if 
there are $\DD$-polynomial maps $f:X \ra Y$ and $g:Y\ra X$
such that $g\circ f = \textup{Id}_X$ and $f\circ g = \textup{Id}_Y$.
\edf

\brm
Let $X=\Spec(K[x_i]_{i\leq q}/(f_j)_{j\in J})$, so that
\[
P_\infty(X) = \Spec(K[d_\albar x_i]_{i\leq q,\albar< \infty}/(f_j)_{j\in J}).
\]
The bijection $\chi$ takes $h \in \Hom_K(\Spec\, K,X)$,
which is determined by $(b_i)_{i\leq q} = h(x_i)_{i\leq q}$
to $H \in \Hom_{(K,\Dun)}(\Spec\, K,P_\infty(X))$ determined by
$(D_\albar b_i)_{i\leq q,\albar < \infty}
=H(d_\albar x_i)_{i\leq q, \albar < \infty}$.
\erm

\bpro
\label{multdpoly}
Let $X$ be a $K$-scheme, and $m<\infty$.  There exists a $\DD$-polynomial
map $\nabla_m:X\ra P_m(X)$ that is a section of the canonical projection
$p_m:P_m(X)\ra X$.  

Let $f:X\ra Y$ be a morphism of $K$-schemes.   Considering $f$ and
$P_m(f)$ as maps on $K$-rational points, the following diagram commutes.
$$
\xymatrix{
P_m(X) \ar[r]^{P_m(f)} & P_m(Y)\\
X\ar[u]^{\nabla_m} \ar[r]^f & Y \ar[u]^{\nabla_m}
}
$$
\epro

\bprf
(Again, the proof is similar to the corresponding result 
for fields with a single derivation.)
By the adjointness of $P_\infty(-)$ and $(-)^!$, there is a natural
bijection
$$
\Hom_K((P_\infty(X))^!,P_m(X))\simeq \Hom_{(K,\Dun)}(P_\infty(X),P_\infty(P_m(X))).
$$ 
Let $f:P_\infty(X)\ra P_\infty(P_m(X))$ be the $\DD$-morphism corresponding
to the canonical projection from $(P_\infty(X))^!$ to $P_m(X)$,
and let $\nabla_m$ be the $\DD$-polynomial map corresponding to $f$.
We show that $\nabla_m$ has the desired properties.

It suffices to check locally, so suppose that $X$ is given as
$\Spec(K[x_i]_{i\leq q}/(f_j)_{j\in J})$.  By Remark~\ref{multafsc}, 
$$
\barr{rcl}
P_m(X) & = & \Spec(K[d_\albar x_i]_{\albar\leq m, i \leq q}/(d_\albar f_j)_{\albar\leq m, j\in J})\\
& & \\
P_\infty(X)& = &\Spec(K[d_\albar x_i]_{\albar < \infty, i \leq q}/(d_\albar f_j)_{\albar < \infty, j\in J})\\
& & \\
P_\infty(P_m(X)) & = & \Spec(K[d_\bebar d_\albar x_i]_{\albar\leq m,\bebar <\infty,i\leq q}/(g_h)_{h\in H})\\
\earr
$$
where $(g_h)_{h\in H}$ is the ideal generated by 
$(d_\bebar d_\albar f_j)_{j\in J, \albar \leq m, l < \infty}$.  The $\DD$-morphism  from 
$P_\infty(X)$ to $P_\infty(P_m(X))$, corresponding to the projection 
morphism from $P_\infty(X)$ to $P_m(X)$ is determined by the $\DD$-algebra
homomorphism
\[
K[d_\bebar d_\albar x_i]_{\albar\leq m,\bebar<\infty,i\leq q}/(g_h)_{h\in H}
\lra
K[d_\albar x_i]_{\albar < \infty, i \leq q}/(d_\albar f_j)_{\albar < \infty, j\in J}
\]
such that
\[ 
 d_\bebar d_\albar x_i 
\mapsto
\tbinom{\al_1+\be_1}{\al_1} \cdots \tbinom{\al_n+\be_n}{\al_n} d_{\albar+\bebar}x_i.
\]
One can then see that this determines the $\DD$-polynomial map from 
$X$ to $P_m(X)$ that takes the closed point $(a_i)_{i\leq n}$ to 
$(D_\albar a_i)_{\albar \leq m, i\leq n}$.  By Remark~\ref{multafsc}, this is a section 
of $\pi_m$.
  
Next we argue that $P_m(f)\circ\nabla_X=\nabla_Y\circ f$.  It 
suffices to prove this for affine schemes, so assume that 
$X=\Spec\, K[\xbar]/I$ and $Y=\Spec\, K[\ybar]/J$.  
Let $S=K[\xbar]/I$ and $R=K[\ybar]/J$, and let $f$ also denote
the homomorphism from $R$ to $S$ corresponding to $f:X\ra Y$.
A $K$-rational point of $X$ corresponds to a homomorphism $g$ from $S$ to $K$, which is
determined by the image of $\xbar$, so we think of a $K$-rational point as
a tuple $\abar= g(\xbar)$ of elements of $K$.  Also, $P_m(X)=\Spec\, \HS^m_{S/(K,\Dun)}$ 
is affine, and $\HS^m_{S/(K,\Dun)}$ is generated 
by $(d_\albar x)_{x\in \xbar, \albar\leq m}$.  We have seen that 
$\nabla_X(\abar )=(\abar, D_1(\abar ), \ldots, D_m(\abar ))$.
To be more precise, $\nabla_X(\abar)$ is the $K$-rational point of $P_m(X)$
that corresponds to the map that sends $d_\albar x\in \HS^m_{S/(K,\Dun)}$ to $D_\albar(g(x))\in K$, 
for each $x\in \xbar, \albar \leq m$.

Let $f(\abar) = \bbar\in Y$, $\bbar = (g\circ f(y))_{y\in\ybar}$.  Again,
$\nabla_Y(\bbar)=(\bbar, D_1(\bbar), \ldots, D_m(\bbar))$.
As a map of $K$-algebras, $P_m(f)$ is the map that sends
$d_\albar y$ to $d_\albar f(y)$, for $y \in \ybar, \albar \leq m$.
Thus,
\[
P_m(f)(\abar, D_1(\abar), \ldots, D_m(\abar))=
(\bbar, D_1(\bbar ), \ldots, D_m(\bbar ))
\]
 as desired.
\eprf

\bibliographystyle{alpha}

\begin{thebibliography}{MvdP03}

\bibitem[Bli05]{Bli03}
Manuel Blickle.
\newblock A short course in geometric motivic integration.
\newblock {\em Preprint}, 2005.
\newblock \texttt{arXiv:math.AG/0507404}.

\bibitem[Bui92]{Bui92}
Alexandru Buium.
\newblock {\em Differential algebraic groups of finite dimension}, volume 1506
  of {\em Lecture Notes in Mathematics}.
\newblock Springer-Verlag, Berlin, 1992.

\bibitem[Bui93]{Bui93}
Alexandru Buium.
\newblock Geometry of differential polynomial functions. {I}. {A}lgebraic
  groups.
\newblock {\em Amer. J. Math.}, 115(6):1385--1444, 1993.

\bibitem[BV95]{BV95}
Alexandru Buium and Jos{\'e}~Felipe Voloch.
\newblock Reduction of the {M}anin map modulo {$p$}.
\newblock {\em J. Reine Angew. Math.}, 460:117--126, 1995.

\bibitem[BV96]{BV96}
Alexandru Buium and Jos{\'e}~Felipe Voloch.
\newblock Lang's conjecture in characteristic {$p$}: an explicit bound.
\newblock {\em Compositio Math.}, 103(1):1--6, 1996.

\bibitem[Cra04]{Cra}
Alastair Craw.
\newblock An introduction to motivic integration.
\newblock In {\em Strings and geometry}, volume~3 of {\em Clay Math. Proc.},
  pages 203--225. Amer. Math. Soc., Providence, RI, 2004.

\bibitem[DL99]{DL}
Jan Denef and Fran{\c{c}}ois Loeser.
\newblock Germs of arcs on singular algebraic varieties and motivic
  integration.
\newblock {\em Invent. Math.}, 135(1):201--232, 1999.

\bibitem[EH00]{EH}
David Eisenbud and Joe Harris.
\newblock {\em The geometry of schemes}, volume 197 of {\em Graduate Texts in
  Mathematics}.
\newblock Springer-Verlag, New York, 2000.

\bibitem[Gil02]{Gil}
Henri Gillet.
\newblock Differential algebra---a scheme theory approach.
\newblock In {\em Differential algebra and related topics (Newark, NJ, 2000)},
  pages 95--123. World Sci. Publishing, River Edge, NJ, 2002.

\bibitem[Har77]{Har}
Robin Hartshorne.
\newblock {\em Algebraic geometry}.
\newblock Springer-Verlag, New York, 1977.
\newblock Graduate Texts in Mathematics, No. 52.

\bibitem[HP00]{HP00}
Ehud Hrushovski and Anand Pillay.
\newblock Effective bounds for the number of transcendental points on
  subvarieties of semi-abelian varieties.
\newblock {\em Amer. J. Math.}, 122(3):439--450, 2000.

\bibitem[Joh85]{Joh}
Joseph Johnson.
\newblock Prolongations of integral domains.
\newblock {\em J. Algebra}, 94(1):173--210, 1985.

\bibitem[Kol73]{Kol73}
E.~R. Kolchin.
\newblock {\em Differential algebra and algebraic groups}.
\newblock Academic Press, New York, 1973.
\newblock Pure and Applied Mathematics, Vol. 54.

\bibitem[Mar06]{Mar96}
David Marker.
\newblock Model theory of differential fields.
\newblock In {\em Model theory of fields}, volume~5 of {\em Lecture Notes in
  Logic}, pages 38--113. Association for Symbolic Logic, La Jolla, CA, second
  edition, 2006.

\bibitem[Mat89]{Mat}
Hideyuki Matsumura.
\newblock {\em Commutative ring theory}, volume~8 of {\em Cambridge Studies in
  Advanced Mathematics}.
\newblock Cambridge University Press, Cambridge, second edition, 1989.
\newblock Translated from the Japanese by M. Reid.

\bibitem[MPS07]{MPS}
Rahim Moosa, Anand Pillay, and Thomas Scanlon.
\newblock Differential arcs and regular types in differential fields.
\newblock {\em \textup{To appear in} J. Reine Angew. Math.}, 2007.

\bibitem[MvdP03]{MvP}
B.~Heinrich Matzat and Marius van~der Put.
\newblock Iterative differential equations and the {A}bhyankar conjecture.
\newblock {\em J. Reine Angew. Math.}, 557:1--52, 2003.

\bibitem[Oku87]{Oku}
K{\^o}taro Okugawa.
\newblock {\em Differential algebra of nonzero characteristic}, volume~16 of
  {\em Lectures in Mathematics}.
\newblock Kinokuniya Company Ltd., Tokyo, 1987.

\bibitem[Pil02]{Pil02}
Anand Pillay.
\newblock Differential fields.
\newblock In {\em Lectures on algebraic model theory}, volume~15 of {\em Fields
  Inst. Monogr.}, pages 1--45. Amer. Math. Soc., Providence, RI, 2002.

\bibitem[Pil04]{Pil04}
Anand Pillay.
\newblock Mordell-{L}ang conjecture for function fields in characteristic zero,
  revisited.
\newblock {\em Compos. Math.}, 140(1):64--68, 2004.

\bibitem[PZ03]{PZ}
Anand Pillay and Martin Ziegler.
\newblock Jet spaces of varieties over differential and difference fields.
\newblock {\em Selecta Math. (N.S.)}, 9(4):579--599, 2003.

\bibitem[Sca97]{Sca97}
Thomas Scanlon.
\newblock The {$abc$} theorem for commutative algebraic groups in
  characteristic {$p$}.
\newblock {\em Internat. Math. Res. Notices}, (18):881--898, 1997.

\bibitem[Sca02]{Sca02}
Thomas Scanlon.
\newblock Model theory and differential algebra.
\newblock In {\em Differential algebra and related topics (Newark, NJ, 2000)},
  pages 125--150. World Sci. Publishing, River Edge, NJ, 2002.

\bibitem[Voj06]{Voj}
Paul Vojta.
\newblock Jets via {H}asse-{S}chmidt derivations.
\newblock In {\em Diophantine Geometry, Proceedings}, pages 335--361. Edizioni
  della Normale, Pisa, 2006.

\bibitem[Zie03]{Zie}
Martin Ziegler.
\newblock Separably closed fields with {H}asse derivations.
\newblock {\em J. Symbolic Logic}, 68(1):311--318, 2003.

\end{thebibliography}

\end{document}